\documentclass[secthm,seceqn,amsthm,ussrhead,10pt]{amsart}
\usepackage[utf8]{inputenc}
\usepackage[english]{babel}

\usepackage{amssymb,amsmath,amsthm,amsfonts,xcolor,enumerate,hyperref,comment,longtable,cleveref}

\usepackage{times}
\usepackage{cite}
\usepackage{pdflscape}
\usepackage{ulem}
\usepackage[mathcal]{euscript}
\usepackage{tikz}
\usepackage{hyperref}
\usepackage{cancel}
\usepackage{stmaryrd}

\usetikzlibrary{arrows}

\newtheorem{theorem}{Theorem}
\newtheorem{lemma}[theorem]{Lemma}

\theoremstyle{definition}
\newtheorem{definition}[theorem]{Definition}

\theoremstyle{remark}
\newtheorem{remark}[theorem]{Remark}

\newenvironment{Proof}[1][Proof.]{\begin{trivlist}
\item[\hskip \labelsep {\bfseries #1}]}{\flushright
$\Box$\end{trivlist}}

\usepackage{stmaryrd}
\usepackage{xcolor}

\newcommand{\aut}[1]{\operatorname{\mathrm{Aut}}{(#1)}}

\newcommand{\A}{\mathbf{A}}

\newcommand{\C}[1]{\mathcal{C}_{#1}}
\newcommand{\R}{\mathcal{R}}

\newcommand{\la}{\langle}
\newcommand{\ra}{\rangle}
\newcommand{\La}{\Big\langle}
\newcommand{\Ra}{\Big\rangle}

\newcommand{\Dl}[2]{[\Delta_{#1#2}]}
\newcommand{\nb}[1]{\nabla_{#1}}

\newcommand{\0}{\theta}
\newcommand{\af}{\alpha}

\newcommand{\gm}{\gamma}

\begin{document}

{\Large\noindent 
Degenerations of nilpotent algebras
\footnote{This work was partially supported by  FAPESP 18/15712-0;
CNPq 451499/2018-2,  404649/2018-1;
RFBR 18-31-20004;
by the President's "Program Support of Young Russian Scientists" (grant MK-2262.2019.1);
by CMUP (UID/MAT/00144/2019), which is funded by FCT with national (MCTES) and European structural funds through the programs FEDER, under the partnership agreement PT2020;
by the Funda\c{c}\~ao para a Ci\^encia e a Tecnologia (Portuguese Foundation for Science and Technology) through the project PTDC/MAT-PUR/31174/2017.}
\footnote{Corresponding author: Ivan Kaygorodov (kaygorodov.ivan@gmail.com)  }
    }

\

{\bf  
        Amir Fern\'andez Ouaridi$^{a}$, 
        Ivan Kaygorodov$^{b,c}$,    
        Mykola Khrypchenko$^{d,e}$
        \&  Yury Volkov$^{f}$\\

    \medskip
}

{\tiny

$^{a}$ Universidad de C\'adiz. Puerto Real, C\'adiz, Spain

$^{b}$ CMCC, Universidade Federal do ABC, Santo Andr\'e, Brazil

$^{c}$ CMUP, Faculdade de Ci\^encias, Universidade do Porto, Porto, Portugal

$^{d}$ Departamento de Matem\'atica, Universidade Federal de Santa Catarina, Florian\'opolis, Brazil

$^{e}$ Departamento de Matem\'atica, Faculdade de Ci\^{e}ncias e Tecnologia, Universidade Nova de Lisboa, Caparica, Portugal

$^{f}$ Saint Petersburg State University, Saint Petersburg, Russia

\

\smallskip

   E-mail addresses:

\smallskip

    Amir Fern\'andez Ouaridi (amir.fernandezouaridi@alum.uca.es)

    Ivan Kaygorodov (kaygorodov.ivan@gmail.com) 
    
    Mykola Khrypchenko (nskhripchenko@gmail.com)    

    Yury Volkov (wolf86\_666@list.ru)

}

\

\noindent{\bf Abstract}: 
{\it 
We give a complete description of degenerations of 
$3$-dimensional nilpotent algebras, 
$4$-dimensional nilpotent commutative algebras
and
$5$-dimensional nilpotent anticommutative algebras over $ \mathbb C$. In particular, we correct several mistakes from the paper ``Contractions of low-dimensional nilpotent Jordan algebras'' by Ancochea Berm\'{u}dez, Fres\'{a}n and Margalef Bentabol.}

\ 

\


\ 


 \ 
 
\section*{Introduction}

There are many results related to the algebraic and geometric 
classification
of low-dimensional algebras in the varieties of Jordan, Lie, Leibniz and 
Zinbiel algebras;
for algebraic classifications  see, for example, \cite{ack, contr11, cfk19, degr3, usefi1, degr2, degr1, demir, gkks, gkk,   ikm19,   ikv18, hac18, kkk18, kpv19, kv16};
for geometric classifications and descriptions of degenerations see, for example, 
\cite{ack, ale, ale2, aleis, maria, contr11, bb14, BC99, cfk19, gkks, gkk19, gkp, GRH, GRH2, ikm19, ikv17, ikv18, kkk18, kpv19, kppv, kpv, kv16,kv17, S90, avdeev, mil, laur03, chouhy, gorb91, gorb93, gorb98, khud15, khud13}.
Here we give the algebraic and geometric classification of  
nilpotent algebras of small dimensions. We also construct the graphs of primary degenerations for the corresponding varieties.

The algebraic classification of nilpotent algebras will be achieved by the calculation of central extensions of algebras from the same variety which have a smaller dimension.
Central extensions of algebras from various varieties were studied, for example, in \cite{ss78,zusmanovich,kkl18,omirov}.
Skjelbred and Sund \cite{ss78} used central extensions of Lie algebras to classify nilpotent Lie algebras.
Using the same method,  
all non-Lie central extensions of  all $4$-dimensional Malcev algebras \cite{hac16},
all non-associative central extensions of all $3$-dimensional Jordan algebras \cite{ha17},
all anticommutative central extensions of $3$-dimensional anticommutative algebras \cite{cfk182},
all central extensions of $2$-dimensional algebras \cite{cfk18}
and some others were described.
One can also look at the classification of
$4$-dimensional nilpotent associative algebras \cite{degr1},
$4$-dimensional nilpotent Novikov algebras \cite{kkk18},
$4$-dimensional nilpotent bicommutative algebras \cite{kpv19},
$5$-dimensional nilpotent restricted Lie agebras \cite{usefi1},
$5$-dimensional nilpotent Jordan algebras \cite{ha16},
$6$-dimensional nilpotent Lie algebras \cite{degr3, degr2},
$6$-dimensional nilpotent Malcev algebras \cite{hac18},
$6$-dimensional nilpotent Tortkara algebras \cite{gkk,gkks},
$6$-dimensional nilpotent binary Lie algebras \cite{ack}.

Degenerations of algebras is an interesting subject, which has been studied in various papers.
In particular, there are many results concerning degenerations of algebras of small dimensions in a  variety defined by a set of identities.
One of important problems in this direction is a description of so-called rigid algebras. 
These algebras are of big interest, since the closures of their orbits under the action of the generalized linear group form irreducible components of the variety under consideration
(with respect to the Zariski topology). 
For example, rigid algebras in the varieties of
all $4$-dimensional Leibniz algebras \cite{ikv17},
all nilpotent $4$-dimensional Novikov algebras \cite{kkk18},
all nilpotent $4$-dimensional bicommutative algebras \cite{kpv19},
all nilpotent $4$-dimensional assosymmetric algebras \cite{ikm19},
all nilpotent $6$-dimensional binary Lie algebras \cite{ack},
and in some other varieties were classified.
There are fewer works in which the full information about degenerations was given for some variety of algebras.
This problem was solved 
for $2$-dimensional pre-Lie algebras in \cite{bb09},  
for $2$-dimensional terminal algebras in \cite{cfk19},
for $3$-dimensional Novikov algebras in \cite{bb14},  
for $3$-dimensional Jordan algebras in \cite{gkp},  
for $3$-dimensional Jordan superalgebras in \cite{maria},
for $3$-dimensional Leibniz and $3$-dimensional anticommutative algebras  in \cite{ikv18},
for $4$-dimensional Lie algebras in \cite{BC99},
for $4$-dimensional Lie superalgebras in \cite{aleis},
for $4$-dimensional Zinbiel and nilpotent $4$-dimensional Leibniz algebras in \cite{kppv},
for nilpotent $5$-dimensional Tortkara algebras in \cite{gkks},
for nilpotent $6$-dimensional Lie algebras in \cite{S90,GRH}, 
for nilpotent $6$-dimensional Malcev algebras in \cite{kpv}, 
for  $2$-step nilpotent $7$-dimensional Lie algebras \cite{ale2}, 
and for all $2$-dimensional algebras in \cite{kv16}.

\section{Preliminaries}

All algebras and vector spaces in this paper are over $\mathbb{C}$ and so we will write simply $\otimes$, ${\rm Hom}$ and $\dim$ instead of $\otimes_{\mathbb{C}}$, ${\rm Hom}_{\mathbb{C}}$ and $\dim_{\mathbb{C}}$.

\subsection{The algebraic classification of nilpotent algebras}\label{algcl}
Let ${\bf A}$ and ${\bf V}$ be an algebra and a vector space and ${\rm Z}^{2}\left( {\bf A},{\bf V}\right)\cong {\rm Hom}({\bf A}\otimes {\bf A},\bf V)$ denote the space of bilinear maps $\theta :{\bf A}\times 
{\bf A}\longrightarrow {\bf V}.$ For $f\in{\rm Hom}({\bf A},{\bf V})$, we introduce $\delta f\in {\rm Z}^{2}\left( {\bf A},{\bf V}\right)$ by the equality $\delta f\left( x,y\right) =f(xy)$ and
define ${\rm B}^{2}\left( {\bf A},{\bf V}\right) =\left\{\delta f \mid f\in {\rm Hom}\left( {\bf A},{\bf V}\right) \right\} $. One
can easily check that ${\rm B}^{2}({\bf A},{\bf V})$ is a linear subspace of ${\rm Z}^{2}\left( {\bf A},{\bf V}\right)$.
Let us define $\rm {H}^{2}\left( {\bf A},{\bf V}\right) $ as the quotient space ${\rm Z}^{2}\left( {\bf A},{\bf V}\right) \big/{\rm B}^{2}\left( {\bf A},{\bf V}\right)$.
The equivalence class of $\theta\in {\rm Z}^{2}\left( {\bf A%
},{\bf V}\right)$ in $\rm {H}^{2}\left( {\bf A},{\bf V}\right)$ is denoted by $\left[ \theta \right]$. We also define $\rm {H}_{\mathcal{C}}^{2}\left( {\bf A},{\bf V}\right)$ as the subspace of $\rm {H}^{2}\left( {\bf A},{\bf V}\right)$ generated by such $\left[ \theta \right]$ that $\theta(x,y)=\theta(y,x)$ for all $x,y\in\bf A$ and  $\rm {H}_{\mathcal{A}}^{2}\left( {\bf A},{\bf V}\right)$  as the subspace of $\rm {H}^{2}\left( {\bf A},{\bf V}\right)$ generated by such $\left[ \theta \right]$ that $\theta(x,x)=0$ for all $x\in\bf A$.

Suppose now that $\dim{\bf A}=m<n$ and $\dim{\bf V}=n-m$. For any
bilinear map $\theta :{\bf A}\times {\bf A}\longrightarrow {\bf V%
}$, one can define on the space ${\bf A}_{\theta }:={\bf A}\oplus 
{\bf V}$ the bilinear product  $\left[ -,-\right] _{%
{\bf A}_{\theta }}$ by the equality $\left[ x+x^{\prime },y+y^{\prime }\right] _{%
{\bf A}_{\theta }}= xy  +\theta \left( x,y\right) $ for  
$x,y\in {\bf A},x^{\prime },y^{\prime }\in {\bf V}$. The algebra ${\bf A}_{\theta }$ is called an $(n-m)$-{\it %
dimensional central extension} of ${\bf A}$ by ${\bf V}$.
It is also clear that${\bf A}_{\theta }$ is nilpotent if and only if  ${\bf A}$ is so.
Moreover, the algebra ${\bf A}_{\theta }$ is (anti)commutative if
and only if ${\bf A}$ is (anti)commutative and $\theta$ is (anti)symmetric.

For a bilinear form $\theta :{\bf A}\times {\bf A}\longrightarrow {\bf V}$, the space $\theta ^{\bot }=\left\{ x\in {\bf A}\mid \theta \left(
x,{\bf A}\right) =\theta \left(
{\bf A},x\right) =0\right\} $ is called the {\it annihilator} of $\theta$.
For an algebra ${\bf A}$, the ideal 
${\rm Ann}\left( {\bf A}\right) =\left\{ x\in {\bf A}\mid x{\bf A} ={\bf A}x =0\right\}$ is called the {\it annihilator} of ${\bf A}$.
One has
\begin{equation*}
{\rm Ann}\left( {\bf A}_{\theta }\right) =\left( \theta ^{\bot }\cap {\rm Ann}\left( 
{\bf A}\right) \right) \oplus {\bf V}.
\end{equation*}
Any $n$-dimensional algebra with non-trivial annihilator can be represented in
the form ${\bf A}_{\theta }$ for some $m$-dimensional algebra ${\bf A}$, an $(n-m)$-dimensional vector space ${\bf V}$ and $\theta \in {\rm Z}^{2}\left( {\bf A},{\bf V}\right)$, where $m<n$ (see \cite[Lemma 5]{hac16}).
Moreover, there is a unique such representation with $m=n-\dim{\rm Ann}({\bf A})$. Note also that the last mentioned equality is equivalent to the condition  $\theta ^{\bot }\cap {\rm Ann}\left( 
{\bf A}\right)=0$. 

Let us pick some $\phi\in {\rm Aut}\left( {\bf A}\right)$, where ${\rm Aut}\left( {\bf A}\right)$ is the automorphism group of  ${\bf A}$.
For $\theta\in {\rm Z}^{2}\left( {\bf A},{\bf V}\right)$, let us define $(\phi \theta) \left( x,y\right) =\theta \left( \phi \left( x\right)
,\phi \left( y\right) \right) $. Then we get an action of ${\rm Aut}\left( {\bf A}\right) $ on ${\rm Z}^{2}\left( {\bf A},{\bf V}\right)$ that induces an action of the same group on $\rm {H}^{2}\left( {\bf A},{\bf V}\right)$. Note that the subspaces $\rm {H}_{\mathcal{C}}^{2}\left( {\bf A},{\bf V}\right)$ and $\rm {H}_{\mathcal{A}}^{2}\left( {\bf A},{\bf V}\right)$ are stable under this action.


\begin{definition}
Let ${\bf A}$ be an algebra and $I$ be a subspace of ${\rm Ann}({\bf A})$. If ${\bf A}={\bf A}_0 \oplus I$
then $I$ is called an {\it annihilator component} of ${\bf A}$.
\end{definition}

For a linear space $\bf U$, the {\it Grassmannian} $G_{s}\left( {\bf U}\right) $ is
the set of all $k$-dimensional linear subspaces of ${\bf U}$. For any $s\ge 1$, the action of ${\rm Aut}\left( {\bf A}\right)$ on $\rm {H}^{2}\left( {\bf A},\mathbb{C}\right)$ induces 
an action of the same group on $G_{s}\left( \rm {H}^{2}\left( {\bf A},\mathbb{C}\right) \right)$.
Let us define
$$
{\bf T}_{s}\left( {\bf A}\right) =\left\{ \bf{W}\in G_{s}\left( \rm {H}^{2}\left( {\bf A},\mathbb{C}\right) \right)\left|\underset{[\theta]\in W}{\cap }\theta^{\bot }\cap {\rm Ann}\left( {\bf A}\right) =0\right.\right\}.
$$
Note that ${\bf T}_{s}\left( {\bf A}\right)$ is stable under the action of ${\rm Aut}\left( {\bf A}\right) $. Note also that $G_{s}\left( \rm {H}_{\mathcal{C}}^{2}\left( {\bf A},\mathbb{C}\right)\right),G_{s}\left( \rm {H}_{\mathcal{A}}^{2}\left( {\bf A},\mathbb{C}\right)\right)\subset G_{s}\left( \rm {H}^{2}\left( {\bf A},\mathbb{C}\right)\right)$.

Let us fix a basis $e_{1},\ldots
,e_{s} $ of ${\bf V}$, and $\theta \in {\rm Z}^{2}\left( {\bf A},{\bf V}\right) $. Then there are unique $\theta _{i}\in {\rm Z}^{2}\left( {\bf A},\mathbb{C}\right)$ ($1\le i\le s$) such that $\theta \left( x,y\right) =\underset{i=1}{\overset{s}{%
\sum }}\theta _{i}\left( x,y\right) e_{i}$ for all $x,y\in{\bf A}$. Note that $\theta ^{\bot
}=\theta^{\bot} _{1}\cap \theta^{\bot} _{2}\cdots \cap \theta^{\bot} _{s}$ in this case.
If   $\theta ^{\bot
}\cap {\rm Ann}\left( {\bf A}\right) =0$, then by \cite[Lemma 13]{hac16} the algebra ${\bf A}_{\theta }$ has a nontrivial
annihilator component if and only if $\left[ \theta _{1}\right] ,\left[
\theta _{2}\right] ,\ldots ,\left[ \theta _{s}\right] $ are linearly
dependent in $\rm {H}^{2}\left( {\bf A},\mathbb{C}\right)$. Thus, if $\theta ^{\bot
}\cap {\rm Ann}\left( {\bf A}\right) =0$ and the annihilator component of ${\bf A}_{\theta }$ is trivial, then $\left\langle \left[ \theta _{1}\right] , \ldots,%
\left[ \theta _{s}\right] \right\rangle$ is an element of ${\bf T}_{s}\left( {\bf A}\right)$.
Now, if $\vartheta\in {\rm Z}^{2}\left( {\bf A},\bf{V}\right)$ is such that $\vartheta ^{\bot
}\cap {\rm Ann}\left( {\bf A}\right) =0$ and the annihilator component of ${\bf A}_{\vartheta }$ is trivial, then by \cite[Lemma 17]{hac16} one has ${\bf A}_{\vartheta }\cong{\bf A}_{\theta }$ if and only if
$\left\langle \left[ \theta _{1}\right] ,\left[ \theta _{2}%
\right] ,\ldots ,\left[ \theta _{s}\right] \right\rangle,
\left\langle \left[ \vartheta _{1}\right] ,\left[ \vartheta _{2}\right] ,\ldots,%
\left[ \vartheta _{s}\right] \right\rangle\in {\bf T}_{s}\left( {\bf A}\right)$ belong to the same orbit under the action of ${\rm Aut}\left( {\bf A}\right) $, where $%
\vartheta \left( x,y\right) =\underset{i=1}{\overset{s}{\sum }}\vartheta
_{i}\left( x,y\right) e_{i}$.

Hence, there is a one-to-one correspondence between the set of $%
{\rm Aut}\left( {\bf A}\right) $-orbits on ${\bf T}_{s}\left( {\bf A}%
\right) $ and the set of isomorphism classes of central extensions of $\bf{A}$ by $\bf{V}$ with $s$-dimensional annihilator and trivial annihilator component.
Consequently to construct all $n$-dimensional central extensions with $s$-dimensional annihilator and trivial annihilator component
of a given $(n-s)$-dimensional algebra ${\bf A}$ one has to describe ${\bf T}_{s}({\bf A})$, ${\rm Aut}({\bf A})$ and the action of ${\rm Aut}({\bf A})$ on ${\bf T}_{s}({\bf A})$ and then
for each orbit under the action of ${\rm Aut}({\bf A})$ on ${\bf T}_{s}({\bf A})$ pick a representative and construct the algebra corresponding to it. If the algebra $\bf{A}$ is (anti)commutative and one wants to construct only (anti)commutative central extensions, then one has to consider ${\bf T}_{s}({\bf A})\cap G_{s}\left( \rm {H}_{\mathcal{C}}^{2}\left( {\bf A},\mathbb{C}\right)\right)$ or ${\bf T}_{s}({\bf A})\cap G_{s}\left( \rm {H}_{\mathcal{A}}^{2}\left( {\bf A},\mathbb{C}\right)\right)$ instead of ${\bf T}_{s}({\bf A})$ correspondingly.

We will use the following auxiliary notation during the construction of central extensions.
Let ${\bf A}$ be an   algebra with  the basis $e_{1},e_{2},\ldots,e_{n}$. In the part devoted to commutative algebras, $\Delta_{ij}:{\bf A}\times {\bf A}\longrightarrow \mathbb{C}$ denotes the symmetric bilinear form  defined by the equalities $\Delta _{ij}\left( e_{i},e_{j}\right)=\Delta _{ij}\left( e_{j},e_{i}\right)=1$
and $\Delta _{ij}\left( e_{l},e_{m}\right) =0$ for $%
\left\{ l,m\right\} \neq \left\{ i,j\right\}$. In this case $\Delta_{ij}$ with $1\leq i \leq j\leq n $ form a basis of the space of symmetric bilinear forms on $\bf{A}$.
 In the part devoted to anticommutative algebras, $\Delta_{ij}:{\bf A}\times {\bf A}\longrightarrow \mathbb{C}$ denotes the antisymmetric bilinear form  defined by the equalities $\Delta _{ij}\left( e_{i},e_{j}\right)=-\Delta _{ij}\left( e_{j},e_{i}\right)=1$
and $\Delta _{ij}\left( e_{l},e_{m}\right) =0$ for $%
\left\{ l,m\right\} \neq \left\{ i,j\right\}$. In this case $\Delta_{ij}$ with $1\leq i < j\leq n $ form a basis of the space of antisymmetric bilinear forms on $\bf{A}$.

\subsection{Degenerations of algebras}
Given an $n$-dimensional vector space ${\bf V}$, the set ${\rm Hom}({\bf V} \otimes {\bf V},{\bf V}) \cong {\bf V}^* \otimes {\bf V}^* \otimes {\bf V}$ 
is a vector space of dimension $n^3$. This space has a structure of the affine variety $\mathbb{C}^{n^3}.$ 
Indeed, let us fix a basis $e_1,\dots,e_n$ of ${\bf V}$. Then any $\mu\in {\rm Hom}({\bf V} \otimes {\bf V},{\bf V})$ is determined by $n^3$ structure constants $c_{i,j}^k\in\mathbb{C}$ such that
$\mu(e_i\otimes e_j)=\sum_{k=1}^nc_{i,j}^ke_k$. A subset of ${\rm Hom}({\bf V} \otimes {\bf V},{\bf V})$ is {\it Zariski-closed} if it can be defined by a set of polynomial equations in the variables $c_{i,j}^k$ ($1\le i,j,k\le n$).

Let $T$ be a set of polynomial identities.
All algebra structures on ${\bf V}$ satisfying polynomial identities from $T$ form a Zariski-closed subset of the variety ${\rm Hom}({\bf V} \otimes {\bf V},{\bf V})$. We denote this subset by $\mathbb{L}(T)$.
The general linear group ${\rm GL}({\bf V})$ acts on $\mathbb{L}(T)$ by conjugation:
$$ (g * \mu )(x\otimes y) = g\mu(g^{-1}x\otimes g^{-1}y)$$ 
for $x,y\in {\bf V}$, $\mu\in \mathbb{L}(T)\subset {\rm Hom}({\bf V} \otimes {\bf V},{\bf V})$ and $g\in {\rm GL}({\bf V})$.
Thus, $\mathbb{L}(T)$ is decomposed into ${\rm GL}({\bf V})$-orbits that correspond to the isomorphism classes of algebras. 
Let $O(\mu)$ denote the ${\rm GL}({\bf V})$-orbit of $\mu\in\mathbb{L}(T)$ and $\overline{O(\mu)}$ its Zariski closure.

Let ${\bf A}$ and ${\bf B}$ be two $n$-dimensional algebras satisfying identities from $T$ and $\mu,\lambda \in \mathbb{L}(T)$ represent ${\bf A}$ and ${\bf B}$ respectively.
We say that ${\bf A}$ {\it degenerates to} ${\bf B}$ and write ${\bf A}\to {\bf B}$ if $\lambda\in\overline{O(\mu)}$.
Note that in this case we have $\overline{O(\lambda)}\subset\overline{O(\mu)}$. Hence, the definition of a degeneration does not depend on the choice of $\mu$ and $\lambda$. If ${\bf A}\not\cong {\bf B}$, then the assertion ${\bf A}\to {\bf B}$ 
is called a {\it proper degeneration}. We write ${\bf A}\not\to {\bf B}$ if $\lambda\not\in\overline{O(\mu)}$.

Let ${\bf A}$ be represented by $\mu\in\mathbb{L}(T)$. Then  ${\bf A}$ is  {\it rigid} in $\mathbb{L}(T)$ if $O(\mu)$ is an open subset of $\mathbb{L}(T)$.
Recall that a subset of a variety is called {\it irreducible} if it cannot be represented as a union of two non-trivial closed subsets. A maximal irreducible closed subset of a variety is called an {\it irreducible component}.
It is well known that any affine variety can be represented as a finite union of its irreducible components in a unique way.
The algebra ${\bf A}$ is rigid in $\mathbb{L}(T)$ if and only if $\overline{O(\mu)}$ is an irreducible component of $\mathbb{L}(T)$.

In the present work we use the methods applied to Lie algebras in \cite{BC99,GRH,GRH2,S90}.
First of all, if ${\bf A}\to {\bf B}$ and ${\bf A}\not\cong {\bf B}$, then $\dim \mathfrak{Der}({\bf A})<\dim \mathfrak{Der}({\bf B})$, where $\mathfrak{Der}({\bf A})$ is the Lie algebra of derivations of ${\bf A}$. We will compute the dimensions of algebras of derivations and will check the assertion ${\bf A}\to {\bf B}$ only for such ${\bf A}$ and ${\bf B}$ that $\dim \mathfrak{Der}({\bf A})<\dim \mathfrak{Der}({\bf B})$. Secondly, if ${\bf A}\to {\bf C}$ and ${\bf C}\to {\bf B}$ then ${\bf A}\to{\bf  B}$. If there is no ${\bf C}$ such that ${\bf A}\to {\bf C}$ and ${\bf C}\to {\bf B}$ are proper degenerations, then the assertion ${\bf A}\to {\bf B}$ is called a {\it primary degeneration}. If $\dim \mathfrak{Der}({\bf A})<\dim \mathfrak{Der}({\bf B})$ and there are no ${\bf C}$ and ${\bf D}$ such that ${\bf C}\to {\bf A}$, ${\bf B}\to {\bf D}$, ${\bf C}\not\to {\bf D}$ and one of the assertions ${\bf C}\to {\bf A}$ and ${\bf B}\to {\bf D}$ is a proper degeneration,  then the assertion ${\bf A} \not\to {\bf B}$ is called a {\it primary non-degeneration}. It suffices to prove only primary degenerations and non-degenerations to describe degenerations in the variety under consideration. It is easy to see that any algebra degenerates to the algebra with zero multiplication. From now on we use this fact without mentioning it.


To prove primary degenerations, we will construct families of matrices parametrized by $t$. Namely, let ${\bf A}$ and ${\bf B}$ be two algebras represented by the structures $\mu$ and $\lambda$ from $\mathbb{L}(T)$ respectively. Let $e_1,\dots, e_n$ be a basis of ${\bf V}$ and $c_{i,j}^k$ ($1\le i,j,k\le n$) be the structure constants of $\lambda$ in this basis. If there exist $a_i^j(t)\in\mathbb{C}$ ($1\le i,j\le n$, $t\in\mathbb{C}^*$) such that $E_i^t=\sum_{j=1}^na_i^j(t)e_j$ ($1\le i\le n$) form a basis of ${\bf V}$ for any $t\in\mathbb{C}^*$, and the structure constants $c_{i,j}^k(t)$ of $\mu$ in the basis $E_1^t,\dots, E_n^t$ satisfy $\lim\limits_{t\to 0}c_{i,j}^k(t)=c_{i,j}^k$, then ${\bf A}\to {\bf B}$. In this case  $E_1^t,\dots, E_n^t$ is called a {\it parametric basis} for ${\bf A}\to {\bf B}$.

To prove primary non-degenerations we will use the following lemma (see \cite{GRH}).

\begin{lemma}\label{main}
Let $\mathcal{B}$ be a Borel subgroup of ${\rm GL}({\bf V})$ and $\mathcal{R}\subset \mathbb{L}(T)$ be a $\mathcal{B}$-stable closed subset.
If ${\bf A} \to {\bf B}$ and ${\bf A}$ can be represented by $\mu\in\mathcal{R}$ then there is $\lambda\in \mathcal{R}$ that represents ${\bf B}$.
\end{lemma}

Each time when we will need to prove some primary non-degeneration $\mu\not\to\lambda$, we will define $\mathcal{R}$ by a set of polynomial equations in structure constants $c_{ij}^k$ in such a way that the structure constants of $\mu$ in the basis $e_1,\dots,e_n$ satisfy these equations. We will omit everywhere the verification of the fact that $\mathcal{R}$ is stable under the action of the subgroup of lower triangular matrices and of the fact that $\lambda\not\in\mathcal{R}$ for any choice of a basis of ${\bf V}.$ 
To simplify our equations, we will use the notation $A_i=\langle e_i,\dots,e_n\rangle,\ i=1,\ldots,n$ and write simply $A_pA_q\subset A_r$ instead of $c_{ij}^k=0$ ($i\geq p$, $j\geq q$, $k\geq r$).

If the number of orbits under the action of ${\rm GL}({\bf V})$ on  $\mathbb{L}(T)$ is finite, then the graph of primary degenerations gives the whole picture. In particular, the description of rigid algebras and irreducible components can be easily obtained. Since the variety of $4$-dimensional nilpotent commutative algebras contains infinitely many non-isomorphic algebras, we have to fulfill some additional work. Let $A(*):=\{A(\alpha)\}_{\alpha\in I}$ be a set of algebras, and let $B$ be another algebra. Suppose that, for $\alpha\in I$, $A(\alpha)$ is represented by the structure $\mu(\alpha)\in\mathbb{L}(T)$ and $B\in\mathbb{L}(T)$ is represented by the structure $\lambda$. Then $A(*)\to B$ means $\lambda\in\overline{\{O(\mu(\alpha))\}_{\alpha\in I}}$, and $A(*)\not\to B$ means $\lambda\not\in\overline{\{O(\mu(\alpha))\}_{\alpha\in I}}$.

Let ${\bf A}(*)$, ${\bf B}$, $\mu(\alpha)$ ($\alpha\in I$) and $\lambda$ be as above. To prove ${\bf A}(*)\to {\bf B}$ it is enough to construct a family of pairs $(f(t), g(t))$ parametrized by $t\in\mathbb{C}^*$, where $f(t)\in I$ and $g(t)\in {\rm GL}({\bf V})$. Namely, let $e_1,\dots, e_n$ be a basis of ${\bf V}$ and $c_{i,j}^k$ ($1\le i,j,k\le n$) be the structure constants of $\lambda$ in this basis. If we construct $a_i^j:\mathbb{C}^*\to \mathbb{C}$ ($1\le i,j\le n$) and $f: \mathbb{C}^* \to I$ such that $E_i^t=\sum_{j=1}^na_i^j(t)e_j$ ($1\le i\le n$) form a basis of ${\bf V}$ for any  $t\in\mathbb{C}^*$, and the structure constants $c_{i,j}^k(t)$ of $\mu\big(f(t)\big)$ in the basis $E_1^t,\dots, E_n^t$ satisfy $\lim\limits_{t\to 0}c_{i,j}^k(t)=c_{i,j}^k$, then ${\bf A}(*)\to {\bf B}$. In this case  $E_1^t,\dots, E_n^t$ and $f(t)$ are called a parametric basis and a {\it parametric index} for ${\bf A}(*)\to {\bf B}$ respectively. In the proofs of this sort, we will put the parametric index in assertion and write $\mu\big(f(t)\big)\to \lambda$ emphasizing that we are proving the assertion $\mu(*)\to\lambda$ using the parametric index $f(t)$.

To prove an assertion of the form ${\bf A}(*)\not\to {\bf B}$, one can use the fact that if $\bf{C}\to \bf{A}(\alpha)$ for any $\alpha\in I$ and $C\not\to \bf{B}$, then ${\bf A}(*)\not\to {\bf B}$.

\section{$3$-dimensional nilpotent algebras}
Thanks to \cite{cfk18}, we have the classification of all $3$-dimensional nilpotent algebras presented in Table A.1 (see Appendix A).
Using this classification and primary degenerations and non-degenerations listed in Tables A.3 and A.4 (see Appendix A) we get the following result.

\begin{theorem}
The variety of $3$-dimensional nilpotent algebras has only one irreducible component defined by 
the rigid algebra $\mathcal N_2.$
The graph of primary degenerations for this variety is given in Figure A.2 (see Appendix A).
\end{theorem}

\section{$4$-dimensional nilpotent commutative algebras}
\subsection{The algebraic classification of $4$-dimensional nilpotent commutative algebras}
Due to the classification of $2$-dimensional algebras, there is only one nontrivial $2$-dimensional nilpotent algebra ${\bf A}_{3}$ (see \cite{kv16}).
The commutative central extensions of ${\bf A}_{3}$ and the trivial $2$-dimensional algebra are described in \cite{cfk18}. In particular, we have the classification of all $3$-dimensional nilpotent commutative algebras, and hence of $4$-dimensional nilpotent commutative algebras with nontrivial annihilator component. This list is formed by the algebras ${\bf A}_{3}, {\bf A}_{08}, {\bf A}_{11}(1)$ and ${\bf A}_{66}$ from \cite{cfk18} that correspond to the algebras $\mathcal{C}_{01}$ -- $\mathcal{C}_{04}$ from Table B.1 of the current paper. Note that the algebra $\mathcal{C}_{04}$ is isomorphic to the algebra ${\bf A}_{66}$ from \cite{cfk18}, but the multiplication table was changed here to simplify our formulas. By the same argument, we have the classification of $4$-dimensional nilpotent commutative algebras with a trivial annihilator component and a $2$-dimensional annihilator. These are the algebras ${\bf A}_{76}(1), {\bf A}_{118}$ and $ {\bf A}_{121}(1)$  from \cite{cfk18} that correspond to the algebras $\mathcal{C}_{05}$ -- $\mathcal{C}_{07}$ from Table B.1 of the current paper.
Moreover, it follows from  \cite{demir} that the only $4$-dimensional central extension of a $3$-dimensional algebra with zero product is the algebra $\mathcal{C}_{08}$ from Table B.1.

The main result of the present subsection is the following theorem.

\begin{theorem}
Let $\bf A$ be a $4$-dimensional nilpotent commutative algebra.
Then $\bf A$ is isomorphic to a unique algebra from the set $\{ \mathcal C_i \}_{i=1, \ldots, 30}$
given in Table B.1 in Appendix B.
\end{theorem}

Due to the paragraph just before the theorem, it is enough to classify $4$-dimensional central extensions of all $3$-dimensional nilpotent commutative algebras with nonzero product.
We will do this in the remaining part of this subsection.

\subsubsection{Automorphism and cohomology groups of of $3$-dimensional nilpotent commutative  algebras}\label{AutHo}

To classify $4$-dimensional central extensions of $3$-dimensional nilpotent commutative algebras we will need the following table describing automorphism groups and cohomology groups of $3$-dimensional nilpotent commutative algebras:

\begin{center}

\begin{tabular}{|c|c|c|c|}
\hline 
$\A$ & $\aut\A$ & ${\rm H}^2_{\mathcal{C}}(\A, \mathbb{C})$\\
\hline

$\mathcal{C}_{01}$ 
& 
$\begin{pmatrix}
x &    0  &  0\\
y &  x^2  &  u\\
z &   0  &  v
\end{pmatrix}$ 
&

$\La 
     \Dl 12, \Dl 13, \Dl 22, \Dl 23,  \Dl 33
\Ra$
\\
\hline

$\mathcal{C}_{02}$

& 
$\begin{pmatrix}
x &    0  &  0\\
0 &  x^2  &  0\\
y &   0  &  x^4
\end{pmatrix}$

&$\La 
     \Dl 12, \Dl 13,  \Dl 23,  \Dl 33
\Ra$\\
\hline

$\mathcal{C}_{03}$ 
& 
$\begin{pmatrix}
x &               0  &  0\\
y &             x^2  &  0\\
z &   2xy  &  x^3
\end{pmatrix}$
&

$\La 
 \Dl 13, \Dl 22, \Dl 23,  \Dl 33
\Ra$\\
\hline

$\mathcal{C}_{04}$ 
& 
$\begin{pmatrix}
x &    0  &  0\\
0 &  v  &  0\\
y &   u  &  xv
\end{pmatrix}$, $\begin{pmatrix}
0 &    x  &  0\\
v &  0  &  0\\
u &   y  &  xv
\end{pmatrix}$ 
&
$\La 
     \Dl 11, \Dl 13, \Dl 22, \Dl 23,  \Dl 33
\Ra$\\

\hline

\hline
\end{tabular}

\end{center}

We give to these algebras the same names as to their $4$-dimensional analogs. The notation used to describe the cohomology is introduced in the end of Subsection \ref{algcl}. The automorphism groups are described by the matrices in the basis $e_1,e_2,e_3$. The variables in these descriptions may take arbitrary values from $\mathbb{C}$ such that the corresponding determinant is nonzero.

\subsubsection{Central extensions of $\mathcal{C}_{01}$}
	Let us use the notation
	$$ 
	\begin{array}{rclrclrclrclrcl}
	\nb 1& = &\Dl 12, & \nb 2& = &\Dl 13, &\nb 3& = &\Dl 23, &\nb 4& = &\Dl 33, &\nb 5& = &\Dl 22.    
	\end{array}
	$$
	Take $\0=\sum_{i=1}^5\af_i\nb i\in {\rm H}^2_{\mathcal{C}}(\mathcal{C}_{01}, \mathbb{C})$.
	If 
	$$
	\phi=
	\begin{pmatrix}
	x &    0  &  0\\
	y &  x^2  &  u\\
	z &   0  &  v
	\end{pmatrix}\in\aut{\mathcal{C}_{01}},
	$$
	then
	$$
	\phi^T
	\begin{pmatrix}
	0          &       \af_1& \af_2\\
	\af_1 &       \af_5& \af_3\\
	\af_2 & \af_3 & \af_4
	\end{pmatrix} 
	\phi=
	\begin{pmatrix}
	\af^{\phi}          &         \af^{\phi}_1& \af^{\phi}_2\\
	\af^{\phi}_1 &         \af^{\phi}_5& \af^{\phi}_3\\
	\af^{\phi}_2 & \af^{\phi}_3 & \af^{\phi}_4	
	\end{pmatrix},
	$$
	i.e. $\phi(\theta)=\sum_{i=1}^5\af^{\phi}_i\nb i$ where
	\begin{align*}
	\af^{\phi}_1 &= x^2(\af_1x + \af_3z + \af_5y),\\
	\af^{\phi}_2 &= \af_1ux + \af_2vx + \af_3vy + \af_4vz + \af_5uy  + \af_3 uz,\\
	\af^{\phi}_3 &= x^2(\af_3v+\af_5u),\\
	\af^{\phi}_4 &= \af_4v^2 + \af_5u^2 +2\af_3uv,\\
	\af^{\phi}_5 &= \af_5x^4.\\
	\end{align*}

	Note that $\mathrm{Ann}({\mathcal C}_{01})=\la e_2, e_3\ra$, and hence the vectors $(\alpha_1,\alpha_5,\alpha_3),(\alpha_1,\alpha_2,\alpha_4)\in\mathbb{C}^3$ have to be linearly independent to give an algebra with a $1$-dimensional annihilator. Let us consider all possible situations.

\begin{enumerate}
\item If $\alpha_5=\alpha_3=0$, then $\alpha_1,\alpha_4\not=0$. Taking $x=\alpha_1^{-\frac{1}{3}}$, $v=\alpha_4^{-\frac{1}{2}}$, $u=y=z=0$, we get $\alpha_1^{\phi}=\alpha_4^{\phi}=1$, and hence we may assume that $\alpha_1=\alpha_4=1$. Then, choosing $x=v=1$, $y=z=0$ and $u=-\alpha_2$, we obtain the subspace $\langle \nabla_1 +\nabla_4 \rangle$.

\item If $\alpha_5=0$, $\alpha_3\not=0$, then we clearly may assume that $\alpha_3=1$. Taking $x=v=1$, $u=-\frac{1}{2}\alpha_4$, $z=-\alpha_1$ and $y=0$, we get $\alpha_1^{\phi}=\alpha_4^{\phi}=0$, and hence we may assume that $\alpha_1=\alpha_4=0$. Then, choosing $x=v=1$, $u=z=0$ and $y=-\alpha_2$, we obtain the subspace $\langle \nabla_3 \rangle$.

\item If $\alpha_5\not=0$, then $\alpha_3^{\phi}=0$, $\alpha_5^{\phi}=1$ for $v=\alpha_5$, $x=\alpha_5^{-\frac{1}{4}}$, $u=-\alpha_3$ and $y=z=0$, and hence we may assume that $\alpha_3=0$ and $\alpha_5=1$. Finally, we have two cases.
\begin{enumerate}
   
   \item If $\alpha_4\not=0$, then we can get $\alpha_4=1$ keeping the equalities $\alpha_3=0$ and $\alpha_5=1$ valid. Then, choosing $x=v=1$, $u=0$, $y=-\alpha_1$ and $z=-\alpha_2$, we obtain the subspace $\langle \nabla_4 +\nabla_5 \rangle$.

   \item If $\alpha_4=0$, then $\alpha_2\not=0$ and we can get $\alpha_2=1$ keeping the equalities $\alpha_3=\alpha_4=0$ and $\alpha_5=1$ valid. Then, choosing $x=v=1$, $u=z=0$ and $y=-\alpha_1$, we obtain the subspace $\langle \nabla_2 +\nabla_5 \rangle$.
\end{enumerate}
\end{enumerate}

It is easy to check that the orbits of obtained subspaces are disjoint. Thus, we get the algebras $ \mathcal{C}_{09}$ -- $\mathcal{C}_{12}.$

\subsubsection{Central extensions of $\mathcal{C}_{02}$}

	Let us use the notation
	$$ 
	\begin{array}{rclrclrclrclrcl}
	\nb 1& = &\Dl 12, & \nb 2& = &\Dl 13, &\nb 3& = &\Dl 23, &\nb 4& = &\Dl 33.    
	\end{array}
	$$
	Take $\0=\sum_{i=1}^4\af_i\nb i\in {\rm H}^2_{\mathcal{C}}(\mathcal{C}_{02}, \mathbb{C})$.
	If 
	$$
	\phi=
	\begin{pmatrix}
	x &    0  &  0\\
	0 &  x^2  &  0\\
    y &   0  &  x^4
	\end{pmatrix}\in\aut{\mathcal{C}_{02}},
	$$
	then
	$$
	\phi^T
	\begin{pmatrix}
	0       &       \af_1   & \af_2\\
	\af_1   &       0       & \af_3\\
	\af_2   &       \af_3   & \af_4
	\end{pmatrix} 
	\phi=
	\begin{pmatrix}
	\af^{\phi}          &         \af^{\phi}_1& \af^{\phi}_2\\
	\af^{\phi}_1 &         0& \af^{\phi}_3\\
	\af^{\phi}_2 & \af^{\phi}_3 & \af^{\phi}_4	
	\end{pmatrix},
	$$
	i.e. $\phi(\theta)=\sum_{i=1}^4\af^{\phi}_i\nb i$ where
	\begin{align*}
	\af^{\phi}_1 &= x^2(\af_1x + \af_3y ),\\
	\af^{\phi}_2 &= x^4(\af_2x+\af_4y),\\
	\af^{\phi}_3 &= \af_3 x^6,\\
	\af^{\phi}_4 &= \af_4x^8.\\
	\end{align*}

	Since $\mathrm{Ann}(\C{02})=\la e_3\ra$, the vector $(\alpha_2,\alpha_3,\alpha_4)\in\mathbb{C}^3$ should be nonzero to give an algebra with a $1$-dimensional annihilator.  Let us consider all possible situations.

\begin{enumerate}
    \item If $\alpha_1=\alpha_3=\alpha_4=0$, then $\alpha_2\not=0$ and we  obtain the subspace $\langle\nabla_2  \rangle$.
    \item If $\alpha_3=\alpha_4=0$ and $\alpha_1\not=0$, then $\alpha_2\not=0$ and choosing $x=\alpha_1^{\frac{1}{2}}\alpha_2^{-\frac{1}{2}}$  we  obtain the subspace $\langle \nabla_1 + \nabla_2  \rangle$.
    \item If $\alpha_1=\alpha_3=0$ and $\alpha_4\not=0$, then choosing $x=\alpha_4$ and $y=-\alpha_2$ we obtain the subspace $\langle\nabla_4  \rangle$.
    \item If $\alpha_3=0$ and $\alpha_1,\alpha_4\not=0$, then choosing $x=\alpha_1^{\frac{1}{5}}\alpha_4^{-\frac{1}{5}}$ and $y=-\alpha_1^{\frac{1}{5}}\alpha_2\alpha_4^{-\frac{6}{5}}$ we obtain the subspace $\langle\nabla_1+\nabla_4  \rangle$.
    \item If $\alpha_3\not=0$, then taking $x=\alpha_3$, $y=-\alpha_1$, we get $\alpha_1^{\phi}=0$, and hence we may assume that $\alpha_1=0$.  Finally we have three cases.
\begin{enumerate}
   \item If $\alpha_2=\alpha_4=0$, then we  obtain the subspace $\langle\nabla_3  \rangle$.
   \item If $\alpha_4=0$ and $\alpha_2\not=0$, then choosing $x=\alpha_1\alpha_2^{-1}$ and $y=0$ we  obtain the subspace $\langle \nabla_2 + \nabla_3  \rangle$.
\item If $\alpha_4\neq0$, then choosing $x=\alpha_3^{\frac{1}{2}}\alpha_4^{-\frac{1}{2}}$ and $y=0$ we obtain the subspace $\langle \alpha \nabla_2 + \nabla_3 + \nabla_4 \rangle$ for some $\alpha\in\mathbb{C}$. Hence, we get the family of subspaces $\langle \alpha \nabla_2 + \nabla_3 + \nabla_4 \rangle$ parameterized by $\alpha$. The subspaces $\langle \alpha \nabla_2 + \nabla_3 + \nabla_4 \rangle$ and $\langle \alpha' \nabla_2 + \nabla_3 + \nabla_4 \rangle$ belong to the same orbit if and only if $\alpha'=-\alpha$.
Thus, the orbits of this subspaces are parameterized by $\alpha\in\mathbb C_{\geq0}$, where $\mathbb{C}_{\geq0}=\{c\in\mathbb{C}\mid {\rm Re}(c)>0\mbox{ or }{\rm Re}(c)=0,{\rm Im}(c)\ge 0\}$.
\end{enumerate}	

\end{enumerate}

It is easy to check that the orbits of obtained subspaces are disjoint. Thus, we get the algebras $\mathcal{C}_{13}$ -- $\mathcal{C}_{18}$ and $\mathcal{C}_{19}(\alpha)$, $\alpha\in \mathbb C_{\geq0}$.

\subsubsection{Central extensions of $\mathcal{C}_{03}$}

	Let us use the notation
	$$ 
	\begin{array}{rclrclrclrclrcl}
	\nb 1& = &\Dl 13, & \nb 2& = &\Dl 22, &\nb 3& = &\Dl 23, &\nb 4& = &\Dl 33.    
	\end{array}
	$$
	Take $\0=\sum_{i=1}^4\af_i\nb i\in {\rm H}^2_{\mathcal{C}}(\mathcal{C}_{03},\mathbb{C})$.
	If 
	$$
	\phi=
	\begin{pmatrix}
	x &    0  &  0\\
	y &  x^2  &  0\\
	z &   2xy  &  x^3
	\end{pmatrix}\in\aut{\mathcal{C}_{03}},
	$$
	then
	$$
	\phi^T
	\begin{pmatrix}
	0       &       0   & \af_1\\
	0   &       \af_2       & \af_3\\
	\af_1   &       \af_3   & \af_4
	\end{pmatrix} 
	\phi=
	\begin{pmatrix}
	\af^{\phi}          &         \beta^{\phi}& \af^{\phi}_1\\
	\beta^{\phi}  &         \af^{\phi}_2 & \af^{\phi}_3\\
	\af^{\phi}_1 & \af^{\phi}_3 & \af^{\phi}_4	
	\end{pmatrix},
	$$
	i.e. $\phi(\theta)=\sum_{i=1}^5\af^{\phi}_i\nb i$ where
	\begin{align*}
	\af^{\phi}_1 &= x^3(\alpha_1x+\alpha_3y+\alpha_4z),\\
	\af^{\phi}_2 &= x^2(\alpha_2x^2+4\alpha_3xy+4\alpha_4y^2),\\
	\af^{\phi}_3 &= x^4(\alpha_3x+2\alpha_4y),\\
	\af^{\phi}_4 &= \af_4x^6.\\
	\end{align*}

Since $\mathrm{Ann}(\C{03})=\la e_3\ra$, the vector $(\alpha_1,\alpha_3,\alpha_4)\in\mathbb{C}^3$ should be nonzero to give an algebra with a $1$-dimensional annihilator.  Let us consider all possible situations.

\begin{enumerate}
    \item If $\alpha_3=\alpha_4=0$, then $\alpha_1\neq 0$ and we  obtain the family of subspace $\langle \nabla_1 + \alpha \nabla_2  \rangle$ parameterized by $\alpha\in\mathbb{C}$. All of them belong to different orbits under the action of $\aut{\mathcal{C}_{03}}$.
    \item If $\alpha_4=0$ and $\alpha_3\not=0$, then taking $x=4\alpha_3$, $y=-\alpha_2$ and $z=0$, we get $\alpha_2^{\phi}=0$, and hence we may assume that $\alpha_2=0$.  Then we have two cases.
\begin{enumerate}
   \item If $\alpha_1=0$, then we  obtain the subspace $\langle\nabla_3  \rangle$.
   \item If $\alpha_1\not=0$, then choosing $x=\alpha_1\alpha_3^{-1}$ and $y=z=0$ we  obtain the subspace $\langle \nabla_1 + \nabla_3  \rangle$.
\end{enumerate}	

\item  If $\alpha_4\neq0$, then  taking $x=2\alpha_4^2$, $y=-\alpha_3\alpha_4$ and $z=\alpha_3^2-2\alpha_1\alpha_4$, we get $\alpha_1^{\phi}=\alpha_3^{\phi}=0$, and hence we may assume that $\alpha_1=\alpha_3=0$. Finally we have two cases.
\begin{enumerate}
   \item If $\alpha_2=0$, then we  obtain the subspace $\langle\nabla_4  \rangle$.
   \item If $\alpha_2\not=0$, then choosing $x=\alpha_2^{\frac{1}{2}}\alpha_4^{-\frac{1}{2}}$ and $y=z=0$ we  obtain the subspace $\langle \nabla_2 + \nabla_4  \rangle$.
\end{enumerate}	
\end{enumerate}

It is easy to check that the orbits of obtained subspaces are disjoint. Thus, we get the algebras $\mathcal{C}_{20}(\alpha)$, $\alpha\in \mathbb C$, and $\mathcal{C}_{21}$ -- $\mathcal{C}_{24}$.

\subsubsection{Central extensions of $\mathcal{C}_{04}$}
	Let us use the notation
	$$ 
	\begin{array}{rclrclrclrclrcl}
	\nb 1& = &\Dl 11, & \nb 2& = &\Dl 13, &\nb 3& = &\Dl 22, &\nb 4& = &\Dl 23, &\nb 5& = &\Dl 33.    
	\end{array}
	$$
	Take $\0=\sum_{i=1}^5\af_i\nb i\in {\rm H}^2_{\mathcal{C}}(\mathcal{C}_{04},\mathbb{C})$.
	If 
	$$
	\phi=
\begin{pmatrix}
x &    0  &  0\\
0 &  v  &  0\\
y &   u  &  xv
\end{pmatrix} \in\aut{\mathcal{C}_{04}},
	$$
	then
	$$
	\phi^T
	\begin{pmatrix}
	\af_1          &       0& \af_2\\
	0 &       \af_3& \af_4\\
	\af_2 & \af_4 & \af_5
	\end{pmatrix} 
	\phi=
	\begin{pmatrix}
	\af_1^{\phi}          &         \af^{\phi}& \af^{\phi}_2\\
	\af^{\phi} &       \af^{\phi}_3& \af^{\phi}_4\\
	\af^{\phi}_2 & \af^{\phi}_4 & \af^{\phi}_5	
	\end{pmatrix},
	$$
	i.e. $\phi(\theta)=\sum_{i=1}^5\af^{\phi}_i\nb i$ where
	\begin{align*}
	\af^{\phi}_1 &=  \alpha_1x^2+2\alpha_2xy+\alpha_5y^2,\\
	\af^{\phi}_2 &=  (\alpha_2x+\alpha_5y)xv,\\
	\af^{\phi}_3 &=   \alpha_3v^2+2\alpha_4vu+\alpha_5u^2,\\
	\af^{\phi}_4 &=  (\alpha_4 v+\alpha_5 u)xv,\\
	\af^{\phi}_5 &= \af_5x^2v^2.\\
	\end{align*}
If 
$$
\phi=
\begin{pmatrix}
0 &    x  &  0\\
v &  0  &  0\\
u &   y  &  xv
\end{pmatrix} \in\aut{\mathcal{C}_{04}},
$$
then $\phi(\theta)=\af^{\phi}_3\nb 1+\af^{\phi}_4\nb 2+\af^{\phi}_1\nb 3+\af^{\phi}_2\nb 4+\af^{\phi}_5\nb 5$, where $\af^{\phi}_i$, $i=1,\dots,5$, are as above.

Since $\mathrm{Ann}(\C{04})=\la e_3\ra$, the vector $(\alpha_2,\alpha_4,\alpha_5)\in\mathbb{C}^3$ should be nonzero to give an algebra with a $1$-dimensional annihilator.  Let us consider all possible situations. 

\begin{enumerate}
    \item If $\alpha_4=\alpha_5=0$, then $\alpha_2\neq0$. Note that the case $\alpha_2=\alpha_5=0$, $\alpha_4\neq0$ can be reduced to this case by interchanging $\alpha_1$ with $\alpha_3$ and $\alpha_2$ with $\alpha_4$.
    Choosing $x=2\alpha_2$, $y=-\alpha_1$, $v=1$ and $u=0$, we get $\alpha_1^{\phi}=0$, and hence we may assume that $\alpha_1=0$. Then we have two cases.
\begin{enumerate}
   \item If $\alpha_3=0$, then we  obtain the subspace $\langle\nabla_2  \rangle$.
   \item If $\alpha_3\not=0$, then choosing $x=1$, $v=\alpha_2\alpha_3^{-1}$ and $y=u=0$ we  obtain the subspace $\langle \nabla_2 + \nabla_3  \rangle$.
\end{enumerate}	
 \item If $\alpha_5=0$ and $\alpha_2,\alpha_4\neq0$, then choosing $v=2\af_2^2\alpha_4$, $x=2\af_2\af_4^2$, $y=-\af_1 \af_4^2$ and $u=-\af_3\af_2^2$ we  obtain the subspace $\langle \nabla_2 + \nabla_4  \rangle$.
\item If $\alpha_5\not=0$, then taking $x=v=\alpha_5$, $y=-\alpha_2$ and $u=-\alpha_4$, we get $\alpha_2^{\phi}=\alpha_4^{\phi}=0$, and hence we may assume that $\alpha_2=\alpha_4=0$. Finally we have three cases.
\begin{enumerate}
   \item If $\alpha_1=\alpha_3=0$, then we  obtain the subspace $\langle\nabla_5  \rangle$.
   \item If $\alpha_3=0$ and $\alpha_1\not=0$, then choosing $x=1$, $v=\alpha_1^{\frac{1}{2}}\alpha_5^{-\frac{1}{2}}$ and $y=u=0$ we  obtain the subspace $\langle \nabla_1 + \nabla_5  \rangle$. Note that the case $\alpha_1=0$, $\alpha_3\neq0$ can be reduced to this case by interchanging $\alpha_1$ with $\alpha_3$ and $\alpha_2$ with $\alpha_4$.
\item If $\alpha_1,\alpha_3\neq0$, then choosing $x=\alpha_3^{\frac{1}{2}}\alpha_5^{-\frac{1}{2}}$, $v=\alpha_1^{\frac{1}{2}}\alpha_5^{-\frac{1}{2}}$ and $y=u=0$ we obtain the subspace $\langle \nabla_1 + \nabla_3 + \nabla_5 \rangle$.
\end{enumerate}	
\end{enumerate}

It is easy to check that the orbits of obtained subspaces are disjoint. Thus, we get the algebras
$\mathcal{C}_{25}$ -- 
$\mathcal{C}_{30}.$

\

\subsection{Degenerations  of $4$-dimensional nilpotent commutative algebras}

\begin{theorem}
The variety of $4$-dimensional nilpotent commutative algebras has only one irreducible component defined by 
the family of algebras $\mathcal C_{19}(\alpha).$
The graph of primary degenerations for this variety is given in Figure B.2 (see Appendix B).
\end{theorem}

\begin{Proof} Tables B.3, B.4 
presented in Appendix B give the proofs for all primary degenerations and non-degenerations.
Table B.5 provides the orbit closures for the families $\C{19}(\af)$ and $\C{20}(\af)$ and, in particular, shows that the whole variety of $4$-dimensional nilpotent commutative algebras coincides with the closure of the set $\C{19}(\af)$ ($\af\in\mathbb{C}_{\ge 0}$).
\end{Proof}

\begin{remark}
    Note that the degenerations $\varphi_5\to\varphi_4$ and $\varphi_4\to\varphi_{10}$ ($\C{09}\to\C{26}$ and $\C{26}\to\C{08}$ in our notation) from~\cite{contr11} are wrong. In fact, we have proved that $\C{09}\not\to\C{26}$ and \ $\C{26}\not\to\C{08}$. Moreover, the authors of~\cite{contr11} affirm that $\varphi_1$ ($\C{20}(1)$ in our notation) cannot degenerate to $\varphi_6$ ($\C{11}$ in our notation), but we have constructed the degeneration $\C{20}(\af)\to \C{11}$ for all $\af\ne 0$.
\end{remark}

\section{$5$-dimensional nilpotent anticommutative  algebras}

\subsection{The algebraic classification of $5$-dimensional nilpotent anticommutative algebras}

\subsubsection{The algebraic classification of $4$-dimensional nilpotent anticommutative algebras}

The classification of $4$-dimensional nilpotent anticommutative algebras (see \cite{cfk182}) is presented in the following table. 
\begin{equation*}
\begin{array}{|l|l|l|l|} 
\hline
\mbox{$\bf A$}  & \mbox{ Multiplication table} & \mbox{${\rm H}_{\mathcal A}^2({\bf A},\mathbb{C})$}  \\ 

\hline
\hline

{\mathcal A}_{01}& e_1e_2 = e_3 &
\la[\Delta_{13}], [\Delta_{14}], [\Delta_{23}], [\Delta_{24}], [\Delta_{34}]  \ra 
  \\
\hline
{\mathcal A}_{02}& e_1e_2=e_3, e_1e_3=e_4 &
\la[\Delta_{14}], [\Delta_{23}], [\Delta_{24}], [\Delta_{34}] \ra \\
\hline

\end{array}
\end{equation*}

In view of \cite{gkks}, all $5$-dimensional anticommutative central extensions of ${\mathcal A}_{01}$ and the algebra with zero multiplication are nilpotent Tortkara algebras that are classified in the same work. Hence, we need to describe only the central extensions of ${\mathcal A}_{02}$.

\subsubsection{$1$-dimensional central extensions of ${\mathcal A}_{02}$}
	Let us use the notation
	$$ 
	\begin{array}{rclrclrclrclrcl}
	\nb 1& = &\Dl 14, & \nb 2& = &\Dl 23, &\nb 3& = &\Dl 24, &\nb 4& = &\Dl 34.    
	\end{array}
	$$
	Take $\0=\sum_{i=1}^4\af_i\nb i\in {\rm H}^2_{\mathcal{A}}(\mathcal{A}_{02},\mathbb{C})$.
	If 
$$
	\phi=
\begin{pmatrix} 
x& 0 & 0 & 0\\
y & z & 0 & 0 \\
u & v & xz  & 0\\
h & g & xv & x^2z
\end{pmatrix}
\in\aut{\mathcal{A}_{02}},
	$$
	then
	$$
	\phi^T\begin{pmatrix} 
0 & 0 & 0 & \alpha_1\\
0 & 0 & \alpha_2 & \alpha_3 \\
0 & -\alpha_2 & 0  & \alpha_4\\
-\alpha_1  & -\alpha_3 & -\alpha_4& 0
\end{pmatrix}
\phi
=
\begin{pmatrix}
0& \alpha^{\phi}& \beta^{\phi}&\alpha_1^{\phi}\\
-\alpha^{\phi}& 0& \alpha_2^{\phi}&\alpha_3^{\phi}\\
-\beta^{\phi}& -\alpha_2^{\phi}& 0& \alpha_4^{\phi}\\
-\alpha_1^{\phi}& -\alpha_3^{\phi}& -\alpha_4^{\phi}& 0
\end{pmatrix},$$
i.e. $\phi(\theta)=\sum_{i=1}^4\af^{\phi}_i\nb i$ where
	$$\begin{array}{lcl}
\alpha_1^{\phi}&=&x^2z(\alpha_1x+\alpha_3y+\alpha_4u), \\
\alpha_2^{\phi}&=&xz(\alpha_2z-\alpha_4g)+vx(\alpha_3z+\alpha_4v), \\
\alpha_3^{\phi}&=&x^2z(\alpha_3z+\alpha_4v),\\
\alpha_4^{\phi}&=&\alpha_4x^3z^2.
\end{array}$$

If $\alpha_4=0$, then the corresponding central extension is a Tortkara algebra due to the results of 
 \cite{gkks}.
If $\alpha_4\neq 0,$ then choosing $x=z=1$, $y=h=0$, $u=-\alpha_1\alpha_4^{-1}, 
v=-\alpha_3\alpha_4^{-1},g=\alpha_2\alpha_4^{-1}$
we get the subspace $\langle \nabla_4 \rangle,$ and hence the algebra
$$\begin{array}{lllll lll}
{\mathcal A}_{11} &:& e_1e_2=e_3,& e_1e_3=e_4,& e_3e_4=e_5. 
\end{array}$$

\subsubsection{The algebraic classification of $5$-dimensional nilpotent  anticommutative algebras}
Since $5$-dimensional nilpotent Tortkara algebras were classified in \cite{gkks}
and we have shown that there is only one $5$-dimensional nilpotent anticommutative non-Tortkara algebra, we get the following theorem.

\begin{theorem}
Let $\bf A$ be a nontrivial $5$-dimensional nilpotent anticommutative algebra. 
Then $\bf A$ is isomorphic to exactly one of the following algebras:

$$\begin{array}{|l|c|lll l|}
 \hline
{\mathcal A}  &  \mathfrak{Der} \ {\mathcal A} & \multicolumn{4}{c|}{\mbox{Multiplication table}}  \\ \hline
{\mathcal A}_{01} & 16 & e_1e_2=e_3. &&&\\ \hline
{\mathcal A}_{02} & 12 & e_1e_2=e_3,& e_1e_3=e_4.&&  \\ \hline
{\mathcal A}_{03} & 15 & e_1e_2=e_4,& e_1e_3=e_5. && \\ \hline
{\mathcal A}_{04} & 10 & e_1e_2=e_3,& e_1e_3=e_4, & e_2e_3=e_5. & \\ \hline
{\mathcal A}_{05} & 11 & e_1e_2=e_5,& e_3e_4=e_5. &&\\ \hline
{\mathcal A}_{06} & 10 & e_1e_2=e_3,& e_1e_4=e_5,& e_2e_3=e_5.& \\ \hline
{\mathcal A}_{07} & 9 & e_1e_2=e_3,& e_3e_4=e_5.&&\\ \hline
{\mathcal A}_{08} & 9 & e_1e_2=e_3,& e_1e_3=e_4,& e_1e_4=e_5.&  \\ \hline
{\mathcal A}_{09} & 8 & e_1e_2=e_3,& e_1e_3=e_4,& e_1e_4=e_5, & e_2e_3=e_5.  \\ \hline
{\mathcal A}_{10} & 7 & e_1e_2=e_3,& e_1e_3=e_4,& e_2e_4=e_5. & \\ \hline
{\mathcal A}_{11} & 6 & e_1e_2=e_3,& e_1e_3=e_4,& e_3e_4=e_5. & \\ \hline
\end{array}$$

\end{theorem}

\subsection{The geometric classification of $5$-dimensional nilpotent  anticommutative algebras}
The degeneration graph for $5$-dimensional nilpotent Tortkara algebras was constructed in \cite{gkks}. In particular, it was shown that this variety has only one irreducible component
 defined by the the  rigid algebra ${\mathcal A}_{10}.$
On the other hand, there is a degeneration ${\mathcal A}_{11} \to {\mathcal A}_{10}$ given by the parametric basis
\[ E_1^t=e_1, \ 
E_2^t=te_2+e_3+\frac{1}{t}e_4, \ 
E_3^t= te_3+e_4, \
E_4^t=te_4, \ 
E_5^t=te_5.\]
This gives us the following theorem. 

\begin{theorem}
The variety of $5$-dimensional nilpotent anticommutative algebras has only one irreducible component defined by the rigid algebra ${\mathcal A}_{11}.$
This variety has the following graph of primary degenerations:

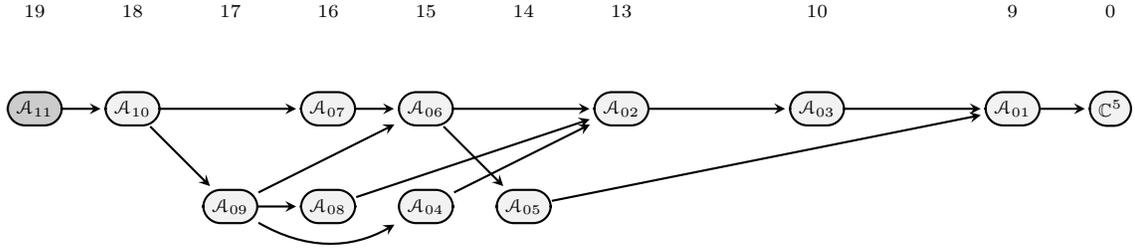
\begin{center}
	
\begin{tikzpicture}[->,>=stealth,shorten >=0.05cm,auto,node distance=1.3cm,
                    thick,main node/.style={rectangle,draw,fill=gray!10,rounded corners=1.5ex,font=\sffamily \scriptsize \bfseries },rigid node/.style={rectangle,draw,fill=black!20,rounded corners=1.5ex,font=\sffamily \scriptsize \bfseries },style={draw,font=\sffamily \scriptsize \bfseries }]

 \node[main node]  (1)                          {$\mathcal A_{09}$ };

  \node[main node] (3) [ right       of=1]       {$\mathcal A_{08}$};

  \node[main node] (2) [ right        of=3]       {$\mathcal A_{04}$};

  \node[main node] (4) [ above        of=2]       {$\mathcal A_{06}$};

  \node[main node] (55) [left          of=4]      {$\mathcal A_{07}$};

  \node[main node] (8) [ right        of=2]       {$\mathcal A_{05}$};

  \node (90) [ above         of=55]      {$16$};
  \node (91) [ left    of=90]      {$17$};
  \node (7) [ left    of=91]      {$18$};
  \node (7xx) [ left    of=7]       {$19$};
  
  \node (92) [ right    of=90]      {$15$};
  \node (93) [ right   of=92]      {$14$};
  \node (94) [ right   of=93]      {$13$};
  \node (94r) [right of=94]{};
  \node (95) [ right   of=94r]      {$10$};
  \node (95r) [right of=95]{};
  \node (96) [ right   of=95r]      {$9$};
  \node (96r) [right of=96]{};
  \node (97) [ right   of=96]      {$0$};

  \node[main node] (t10) [ below         of=7]       {$\mathcal A_{10}$};

  \node[rigid node] (t11)    [ left         of=t10]                       {$\mathcal A_{11}$ };

  \node[main node] (5) [ below         of=94]       {$\mathcal A_{02}$};

\node (5r) [right of=5]{};

  \node[main node] (6) [ right       of=5r]       {$\mathcal A_{03}$};
  
\node (6r) [right of=6]{};

  \node[main node] (7) [ right        of=6r]       {$\mathcal A_{01}$};
  
\node (7r) [right of=7]{};

  \node[main node] (111) [ right   of=7]       {$\mathbb C^5$};
 
  \path[every node/.style={font=\sffamily\small}]
    (1) edge  [bend right] node[left] {} (2)
    (1) edge   node[left] {} (3)
    (1) edge   node[left] {} (4)

    (2) edge   node[left] {} (5)
    
    (3) edge   node[left] {} (5)
    
    (4) edge   node[left] {} (8)
    (4) edge     node[left] {} (5)
    
    (5) edge   node[left] {} (6)

    (6) edge   node[left] {} (7)

    (8) edge  node[left] {} (7)
    
    (7) edge   node[left] {} (111)
    (t10) edge   node[left] {} (1)
    (t10) edge   node[left] {} (55)
    (t11) edge   node[left] {} (t10)

    (55) edge   node[left] {} (4);

\end{tikzpicture}

\end{center}

\end{theorem}

\section*{Appendix A. $3$-dimensional nilpotent algebras}

$$\begin{array}{|l|c|lll|}
\hline

\multicolumn{5}{|c|}{\textrm{{\bf Table A.1. $3$-dimensional nilpotent  algebras}}}  \\

\hline
\A  & \mathfrak{Der} \ \A  & \multicolumn{3}{c|}{\mbox{Multiplication table}} \\ \hline

\mathcal{N}_1 &  2 & e_1 e_1 = e_2, & e_2 e_2=e_3 & \\ \hline
\mathcal{N}_2 &  1 & e_1 e_1 = e_2, & e_2 e_1= e_3, & e_2 e_2=e_3  \\ \hline
\mathcal{N}_3 &  3 & e_1 e_1 = e_2, & e_2 e_1=e_3 &  \\\hline
\mathcal{N}_4(\alpha) &  3 & e_1 e_1 = e_2, & e_1 e_2=e_3, & e_2 e_1=\alpha e_3  \\ \hline
\mathcal{N}_5 &  5  & e_1 e_1 = e_2 && \\ \hline
\mathcal{N}_6 &  4  & e_1 e_1 = e_3, & e_2 e_2=e_3  &\\ \hline
\mathcal{N}_7 &  6  & e_1 e_2=e_3, & e_2 e_1=-e_3  &\\ \hline
\mathcal{N}_8(\alpha) & 4 & e_1 e_1 = \alpha e_3, & e_2 e_1=e_3,  & e_2 e_2=e_3 \\ \hline

\end{array}$$

\ 

\begin{center}
		
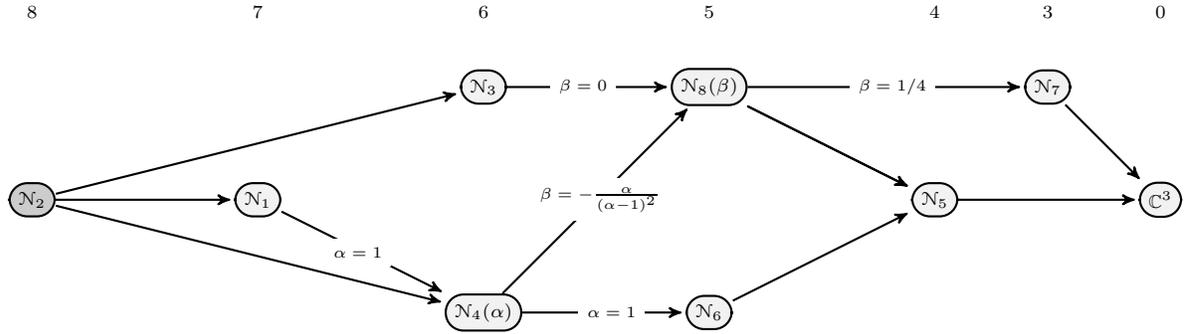
		{\bf Figure A.2. The graph of degenerations of $3$-dimensional nilpotent algebras}
		$$\begin{tikzpicture}[->,>=stealth',shorten >=0.05cm,auto,node distance=1.5cm,
		thick,main node/.style={rectangle,draw,fill=gray!10,rounded corners=1.5ex,font=\sffamily \scriptsize \bfseries },rigid node/.style={rectangle,draw,fill=black!20,rounded corners=1.5ex,font=\sffamily \scriptsize \bfseries },style={draw,font=\sffamily \scriptsize \bfseries }]
		\node (1)       {$8$};
		\node            (r1)  [below of =1,yshift=2cm]   {};                     
		\node            (rr1)  [below of =r1]   {};  
		\node            (rrr1)  [below of =rr1]   {}; 
		\node[rigid node] (N2)  [below of =rr1]                       {$\mathcal{N}_2$ };
		
		\node (2x) [right  of=1]       {};
		\node (2) [right  of=2x]       {$7$};
		\node            (r2)  [below of =2,yshift=2cm]   {};                     
		\node            (rr2)  [below of =r2]   {};  
		\node            (rrr2)  [below of =rr2]   {};  
		\node[main node] (N1)  [below of =rr2]                       {$\mathcal{N}_1$ };
		
		\node (3x) [right  of=2]       { };
		\node (3) [right  of=3x]       {$6$};
		
		\node            (r3)  [below of =3,yshift=2cm]   {};                     
		\node[main node] (N3)  [below of =r3]                       {$\mathcal{N}_3$ };
		
		\node            (rr3)  [below of =N3]   {};                     
		\node            (rrr3)  [below of =rr3]   {}; 
		\node            (rrrr3)  [below of =rrr3]   {}; 
		\node[main node] (N4)  [below of =rr3]                       {$\mathcal{N}_4(\alpha)$ };
		
		\node (3a) [right  of=3]      {};
		
		\node (4) [right  of=3a]      {$5$};
		
		\node            (r4)  [below of =4,yshift=2cm]   {};    
		\node[main node] (N8)  [below of =r4]                       {$\mathcal{N}_{8}(\beta)$};
		
		\node            (rN8)  [below of =N8]   {};		
		
		\node            (rrN8)  [below of =rN8]   {};				
		\node            (rrrN8)  [below of =rrN8]   {};				
		\node[main node] (N6)  [below of =rN8]                       {$\mathcal{N}_6$ };		
		
		\node (5x) [right  of=4]      {};
		\node (5) [right  of=5x]      {$4$};

		\node            (r5)  [below of =5,yshift=2cm]   {};   				
		\node            (rr5)  [below of =r5]   {}; 
		\node            (rrr5)  [below of =rr5]   {}; 
		\node            (rrrr5)  [below of =rrr5]   {}; 
		\node            (rrrrr5)  [below of =rrrr5]   {}; 
		
		\node[main node] (N5)  [below of =rr5]                       {$\mathcal{N}_5$ };
		
		\node (6) [right  of=5]      {$3$};
		
		\node            (r6)  [below of =6,yshift=2cm]   {};  
		\node            (rr6)  [below of =r6]   {};  
		\node[main node] (N7)  [below of =r6]                       {$\mathcal{N}_7$ };
		
		\node (9x) [right  of=6]      {};
		\node (9) [right  of=6]      {$0$};
		
		\node            (r9)  [below of =9,yshift=2cm]   {};  
		\node            (rr9)  [below of =r9]   {};  
		\node            (rrr9)  [below of =rr9]   {};  
		
		\node[main node] (C3)  [below of =rr9]                       {$\mathbb C^3$ };
		
		\path[every node/.style={font=\sffamily\small}]

		(N2)   edge  (N1) 
		(N2)   edge  (N3)
		(N2)   edge  (N4)

		(N6)   edge  (N5)
		(N8)   edge  (N5)
		(N8)   edge  (N5)

		(N1)   edge [bend right=0] node[above=0, right=-15, fill=white]{\tiny $\alpha=1$ }    node {} (N4) 
		
		(N3)   edge [bend right=0] node[above=0, right=-15, fill=white]{\tiny $\beta=0$ }    node {} (N8) 
		
		(N4)   edge [bend right=0] node[above=0, right=-25, fill=white]{\tiny $\beta=-\frac{\alpha}{(\alpha-1)^2}$ }    node {} (N8)    
		
		(N8)   edge [bend right=0] node[above=0, right=-14, fill=white]{\tiny $\beta=1/4$ }    node {} (N7)  
		(N4)   edge [bend right=0] node[above=0, right=-10, fill=white]{\tiny $\alpha=1$ }    node {} (N6)

		(N5)   edge  (C3)
		(N7)   edge  (C3);

		\end{tikzpicture}
		$$
	\end{center}

$$\begin{array}{|rcl|lll|}
\hline

\multicolumn{6}{|c|}{\textrm{{\bf Table A.3. Degenerations of $3$-dimensional nilpotent  algebras}}}  \\
\hline
\multicolumn{3}{|c|}{\textrm{Degeneration}}  & \multicolumn{3}{|c|}{\textrm{Parametric basis}} \\
\hline
\hline
\mathcal{N}_1 &\to& \mathcal{N}_4(1)  & 
E_{1}^{t}=e_1+\frac{1}{t}e_2,& 
E_{2}^{t}=e_2 + \frac{1}{t^2}e_3, & 
E_{3}^{t}=\frac{1}{t}e_3 \\


\hline
\mathcal{N}_2 &\to& \mathcal{N}_1& 
E_{1}^{t}= \frac{1}{t}e_1, & 
E_{2}^{t}=\frac{1}{t^2}e_2, &
E_{3}^{t}=\frac{1}{t^4}e_3 \\

\hline
\mathcal{N}_2 &\to& \mathcal{N}_3 & 
E_{1}^{t}= te_1, &
E_{2}^{t}=t^2 e_2, &
E_{3}^{t}=t^3e_3 \\

\hline
\mathcal{N}_2 &\to& \mathcal{N}_4(\alpha\neq 1) & 
E_1^t=te_1+\frac{t}{\alpha-1}e_2,&
E_2^t=t^2e_2+\frac{\alpha}{(\alpha-1)^2}t^2e_3,&
E_3^t=\frac{t^3}{\alpha-1}e_3,\\



\hline
\mathcal{N}_3 &\to& \mathcal{N}_8(0) & 
E_{1}^{t}=te_1, &
E_{2}^{t}=te_1+e_2, &
E_{3}^{t}=te_3 \\

\hline
\mathcal{N}_4(1) &\to& \mathcal{N}_6 & 

E_1^t=te_1, &
E_2^t=ite_1-ite_2+ite_3,&
E_3^t=t^2e_3,\\

\hline

\mathcal{N}_4(\alpha\neq 0,1) &\to& \mathcal{N}_8 \left(-\frac{\alpha}{(\alpha-1)^2}\right)& 
E_1^t=te_1,&
E_2^t= \frac{\alpha-1}{\alpha}te_1+\frac{t}{\alpha-1}e_2+\frac{\alpha t}{(\alpha-1)^3}e_3,&
E_3^t=t^2e_3,\\

\hline
\mathcal{N}_4(0) &\to& \mathcal{N}_8(0) & 
E_1^t=e_2,&
E_2^t= te_1+e_2,&
E_3^t=te_3,\\
\hline


\mathcal{N}_6  &\to& \mathcal{N}_5 & 
E_{1}^{t}= e_2, & 
E_{2}^{t}= e_3, &
E_{3}^{t}=te_1 \\

\hline

\mathcal{N}_8(\alpha) &\to& \mathcal{N}_5 & 
E_{1}^{t}= e_2, & 
E_{2}^{t}= e_3, &
E_{3}^{t}=te_1 \\
\hline

\mathcal{N}_8(1/4) &\to& \mathcal{N}_7 & 
E_{1}^{t}= t e_2, & 
E_{2}^{t}=2 e_1 - e_2, &
E_{3}^{t}=te_3 \\

\hline


\end{array}$$

$$\begin{array}{|rcl|l|}
\hline
\multicolumn{4}{|c|}{\textrm{{\bf Table A.4. Non-degenerations of $3$-dimensional nilpotent  algebras}}}  \\
\hline
\hline
\multicolumn{3}{|c|}{\textrm{Non-degeneration}} & \multicolumn{1}{|c|}{\textrm{Arguments}}\\
\hline
\hline

\hline
\mathcal{N}_1  &\not \to& \mathcal{N}_7, \mathcal{N}_8(\beta) & 
\mbox{$\mathcal{N}_1$\rm  is commutative while $\mathcal{N}_7$ and $\mathcal{N}_8(\beta)$ are not commutative}\\
 \hline

\mathcal{N}_3  &\not \to& \mathcal{N}_6,\mathcal{N}_7,\mathcal{N}_8(\alpha\neq0) & {\mathcal R}=\left\{  
A_1A_2=0\right\} \\

\hline
\mathcal{N}_4(\alpha) &\not \to& 

\begin{array}{ll}
\mathcal{N}_6 &(\mbox{for \ } \alpha\neq 1), \\
\mathcal{N}_8(\beta) &(\mbox{for \ } \beta \neq- \frac{\alpha}{(\alpha-1)^2}),\\
\mathcal{N}_7 & (\mbox{for \ } \alpha \neq -1)
\end{array}

& {\mathcal R}=\left\{    
\begin{array}{l} A_3A_1+A_1A_3+A_2A_2=0,\\
A_1A_2+A_2A_1\subset A_3,c_{21}^3=\alpha c_{12}^3
\end{array}
\right\} \\

\hline




\hline
\mathcal{N}_8(\alpha\neq 1/4) &\not \to& \mathcal{N}_7 & {\mathcal R}=
\left\{
\begin{array}{l} A_3A_1+A_1A_3=0, A_1A_1\subseteq A_3\\ c_{11}^3c_{22}^3=\alpha (c_{21}^3-c_{12}^3)^2+c_{12}^3c_{21}^3
\end{array}
\right\}
 \\

\hline 
\end{array}$$

\section*{Appendix B. $4$-dimensional nilpotent commutative algebras}

\begin{center}

$$	\begin{array}{|l|c|lllll|}
			\hline 
			\multicolumn{7}{|c|}{\textrm{{\bf Table B.1. $4$-dimensional nilpotent commutative  algebras}}}  \\
			\hline
			\hline 

			\A & \mathfrak{Der} \ \A & 			\multicolumn{5}{|c|}{\textrm{Multiplication table}}  \\

			\hline
			
			\mathcal{C}_{01} & 10 & e_1 e_1 = e_2 &&&&\\
			
			\hline
			
			\mathcal{C}_{02} & 5 & e_1 e_1 = e_2, & e_2 e_2=e_3 &&&\\
			
			\hline
			
			\mathcal{C}_{03} & 6 & e_1 e_1 = e_2, & e_1 e_2=e_3&&&
			\\
			\hline
			
			\mathcal{C}_{04} & 8 & e_1 e_2=e_3&&&& 
			\\
			\hline
			
			\mathcal{C}_{05} & 4 & e_1 e_1 = e_2, & e_1 e_2=e_4, & e_2 e_2=e_3&&
			\\
			\hline
			
			\mathcal{C}_{06} & 6 & e_1 e_1 = e_3, & e_2 e_2=e_4 &&&
			\\
			\hline
			
			\mathcal{C}_{07} & 7 & e_1 e_1 = e_3, & e_1 e_2=e_4 &&&
			\\
			\hline
			
			\mathcal{C}_{08} & 7 & e_1 e_1 = e_4, & e_2 e_3=e_4 &&&
			\\
			\hline

			\mathcal{C}_{09} & 4 & e_1 e_1 = e_2,& e_2e_3=e_4 &&&
			\\
			\hline
			
			\mathcal{C}_{10} & 4 & e_1 e_1 = e_2,& e_1e_3=e_4,& e_2e_2=e_4&&
			\\
			\hline
			
			\mathcal{C}_{11} & 5 & e_1 e_1 = e_2,& e_1e_2=e_4,&e_3e_3=e_4 &&
			\\
			\hline
			
			\mathcal{C}_{12} & 3 & e_1 e_1 = e_2,& e_2e_2=e_4,& e_3e_3=e_4 &&
			\\
			\hline
			
			\mathcal{C}_{13} & 3 & e_1 e_1 = e_2,& e_1e_3=e_4, & e_2 e_2=e_3&& 
			\\
			\hline
			
			\mathcal{C}_{14} & 2 & e_1 e_1 = e_2,&e_1e_2= e_4, & e_1e_3=e_4, & e_2 e_2=e_3&
			\\
			\hline
			
			\mathcal{C}_{15} & 2 & e_1 e_1 = e_2, & e_2 e_2=e_3,  &e_2e_3=e_4 &&
			\\
			\hline
			
			\mathcal{C}_{16} & 1 & e_1 e_1 = e_2,& e_1e_3= e_4, & e_2 e_2=e_3, &e_2e_3=e_4 & 
			\\
			\hline
			
			\mathcal{C}_{17} & 2 & e_1 e_1 = e_2, & e_2 e_2=e_3, &e_3e_3=e_4  &&
			\\
			\hline
			
			\mathcal{C}_{18} & 1 & e_1 e_1 = e_2,& e_1e_2= e_4, & e_2 e_2=e_3, &e_3e_3=e_4 &
			\\
			\hline
			
			\mathcal{C}_{19}(\alpha \in  \mathbb{C}_{\geq 0}) & 1 & e_1 e_1 = e_2,& e_1e_3=\alpha e_4, & e_2 e_2=e_3, &e_2e_3= e_4, & e_3e_3=e_4 

			\\
			\hline
			
			\mathcal{C}_{20}(\alpha) & 4 & e_1 e_1 = e_2, & e_1 e_2=e_3,& e_1e_3=e_4, & e_2e_2=\alpha e_4 &
			\\
			\hline
			
			\mathcal{C}_{21} & 3 & e_1 e_1 = e_2, & e_1 e_2=e_3,& e_2e_3= e_4 &&
			\\
			\hline
			
			\mathcal{C}_{22} & 2 & e_1 e_1 = e_2, & e_1 e_2=e_3,& e_1e_3= e_4, & e_2e_3= e_4 &
			\\
			\hline
			
			\mathcal{C}_{23} & 2 & e_1 e_1 = e_2, & e_1 e_2=e_3,&  e_3e_3= e_4 &&
			\\
			\hline
			
			\mathcal{C}_{24} & 1 & e_1 e_1 = e_2, & e_1 e_2=e_3,& e_2e_2= e_4, & e_3e_3= e_4 &
			\\
			\hline
			
			\mathcal{C}_{25} & 3 & e_1 e_2 = e_3, & e_1 e_3=e_4,& e_2e_3=e_4 && 
			\\
			\hline
			
			\mathcal{C}_{26} & 5 & e_1 e_2 = e_3, & e_1 e_3=e_4&&&
			\\
			\hline
			
			\mathcal{C}_{27} & 4 & e_1 e_2 = e_3,& e_1e_3=e_4, & e_2 e_2=e_4&&
			\\
			\hline
			
			\mathcal{C}_{28} & 4  & e_1 e_2 = e_3, & e_3 e_3=e_4&&&
			\\
			\hline
						
			\mathcal{C}_{29} & 2 & e_1 e_1 = e_4, & e_1 e_2=e_3,& e_2e_2=e_4,& e_3e_3=e_4 &
			\\
			\hline
						
			\mathcal{C}_{30} & 3 & e_1 e_1 = e_4, & e_1 e_2=e_3,& e_3e_3=e_4&&
			\\
			\hline
		\end{array}$$ 

\end{center}

\newpage
\begin{center}
	{\bf Figure B.2.  The graph of degenerations of $4$-dimensional nilpotent commutative algebras}
	
	\
	
	\begin{tikzpicture}[->,>=stealth,shorten >=0.05cm,auto,node distance=1.3cm,
	thick,main node/.style={rectangle,draw,fill=gray!10,rounded corners=1.5ex,font=\sffamily \scriptsize \bfseries },rigid node/.style={rectangle,draw,fill=black!20,rounded corners=1.5ex,font=\sffamily \scriptsize \bfseries },style={draw,font=\sffamily \scriptsize \bfseries }]

	\node (15) at (0,24) {$15$};
	\node (14) at (0,21) {$14$};
	\node (13) at (0,18) {$13$};
	\node (12) at (0,15) {$12$};
	\node (11) at (0,12) {$11$};
	\node (10) at (0,9) {$10$};
	\node (9)  at (0,7) {$9$};
	\node (8)  at (0,5) {$8$};
	\node (6)  at (0,3) {$6$};
	\node (0)  at (0,1) {$0$};

	\node[main node] (c16) at (-4,24) {$\C{16}$};
	\node[main node] (c18) at (-7,24) {$\C{18}$};
	\node[rigid node] (c19) at (-10,24) {$\C{19}(\alpha)$};
	\node[main node] (c24) at (-13,24) {$\C{24}$};
	
	\node[main node] (c14) at (-1,21) {$\C{14}$};
	\node[main node] (c15) at (-7,21) {$\C{15}$};
	\node[main node] (c17) at (-4,21) {$\C{17}$};
	\node[main node] (c22) at (-10,21) {$\C{22}$};
	\node[main node] (c23) at (-13,21) {$\C{23}$};
	\node[main node] (c29) at (-16,21) {$\C{29}$};
	
	\node[main node] (c21) at (-14.5,18) {$\C{21}$};
	\node[main node] (c25) at (-8.5,18) {$\C{25}$};
    \node[main node] (c30) at (-5.5,18) {$\C{30}$};
	\node[main node] (c13) at (-2.5,18) {$\C{13}$};
	\node[main node] (c12) at (-11.5,18) {$\C{12}$};

	\node[main node] (c05) at (-1,15) {$\C{05}$};
	\node[main node] (c28) at (-4,15) {$\C{28}$};
	\node[main node] (c10) at (-7,15) {$\C{10}$};
	\node[main node] (c20) at (-10,15) {$\C{20}(\alpha)$};
	\node[main node] (c27) at (-13,15) {$\C{27}$};
	\node[main node] (c09) at (-16,15) {$\C{09}$};

	\node[main node] (c02) at (-2.5,12) {$\C{02}$};
	\node[main node] (c11) at (-8.5,12) {$\C{11}$};
	\node[main node] (c26) at (-14.5,12) {$\C{26}$};
	
	\node[main node] (c03) at (-11.5,9) {$\C{03}$};
	\node[main node] (c06) at (-5.5,9) {$\C{06}$};
	
	\node[main node] (c07) at (-5.5,7) {$\C{07}$};
	\node[main node] (c08) at (-11.5,7) {$\C{08}$};
	
	\node[main node] (c04) at (-8.5,5) {$\C{04}$};
	
	\node[main node] (c01) at (-8.5,3) {$\C{01}$};
	
	\node[main node] (CC)  at (-8.5,1) {$\mathbb{C}^4$};
	
	\path[every node/.style={font=\sffamily\small}]
    (c01) edge   node[left] {} (CC)
    (c02) edge   node[left] {} (c03)
    (c02) edge   node[left] {} (c06)
    (c03) edge   node[left] {} (c07)
    (c04) edge   node[left] {} (c01)
    (c05) edge   node[left] {} (c02)
    (c06) edge   node[left] {} (c07)
    (c07) edge   node[left] {} (c04)
    (c08) edge   node[left] {} (c04)
    (c09) edge   [bend left=12] node[left] {} (c11)
    (c10) edge   node[left] {} (c02)
    (c10) edge   node[left] {} (c11)
    (c11) edge   node[left] {} (c08)
    (c11) edge   node[left] {} (c03)
    (c11) edge   node[left] {} (c06)
    (c12) edge   node[left] {} (c10)
    (c12) edge   node[left] {} (c09)
    (c13) edge   [bend left=5] node[left] {} (c09)
    (c13) edge   node[left] {} (c10)
    (c13) edge   node[above=-10, right=-45, fill=white]{\tiny $\af=-1$ }    node {}  (c20)
    (c14) edge   node[left] {} (c05)
    (c14) edge   node[left] {} (c13)
    (c14) edge   node[left] {} (c25)
    (c14) edge   [bend left=-12] node[left] {} (c20)
    (c15) edge   node[left] {} (c05)
    (c15) edge   node[left] {} (c12)
    (c15) edge   [bend left=-25] node[left] {} (c20)
    (c15) edge   node[left] {} (c25)
    (c16) edge   node[left] {} (c14)
    (c16) edge   node[left] {} (c15)
    (c16) edge   node[left] {} (c22)
    
    (c17) edge   [bend left=-18] node[left] {} (c05)
    (c17) edge   node[left] {} (c12)
    (c17) edge   node[left] {} (c13)
    (c17) edge   [bend right=0]    node {} (c20)
    (c17) edge     node {} (c25)
    
    (c18) edge   node[left] {} (c14)
    (c18) edge   node[left] {} (c17)
    (c18) edge   node[left] {} (c22)
   
   (c19) edge   node[left] {} (c14)
   (c19) edge   node[above=12, right=-23, fill=white]{\tiny $\alpha= 0$ } node[left] {} (c15)
   (c19) edge   node[left] {} (c17)
   (c19) edge   node[above=12, right=-10, fill=white]{\tiny $\alpha\not= 0$ } node[left] {} (c22)
   (c19) edge   node[above=12, right=-0, fill=white]{\tiny $\alpha= 0$ } node[left] {} (c23)
   (c19) edge   node[above=25, right=30, fill=white]{\tiny $\alpha= 0$ } node[left] {} (c29)

    (c20) edge   node[above=0, right=-10, fill=white]{\tiny $\alpha\neq 0$ }    node {}  (c11)
    (c20) edge   node[above=0, right=-10, fill=white]{\tiny $\alpha= 0$ }    node {}  (c26)
    (c21) edge   node[left] {} (c27)
    (c21) edge   [bend left=0] node[above=-7, right=-0, fill=white]{\tiny $\alpha= 4$ }    node {}  (c20)
    (c21) edge   node[left] {} (c09)
    (c22) edge   [bend left=12] node[left] {} (c05)
    (c22) edge   node[left] {} (c12)
    (c22) edge     node {}  (c20)
    (c22) edge   node[left] {} (c21)
    (c22) edge   node[left] {} (c25)
    (c23) edge   [bend left=17] node[left] {} (c05)
    (c23) edge   [bend left=25] node[left] {} (c20)
    (c23) edge [bend left=0]  node[left] {} (c30)
    (c24) edge   node[left] {} (c22)
    (c24) edge   node[left] {} (c23)
    (c24) edge   node[left] {} (c29)
    (c25) edge   node[left] {} (c09)
    (c25) edge   [bend left=-0] node[left] {} (c27)
    (c26) edge   node[left] {} (c03)
    (c27) edge   node[left] {} (c26)
    (c27) edge   node[left] {} (c11)
    (c28) edge   node[left] {} (c02)
    (c28) edge   node[left] {} (c11)
    (c29) edge   node[left] {} (c12)
    (c29) edge   node[left] {} (c25)
   (c29) edge   node[left] {} (c30)
 
    (c30) edge  [bend left=0] node[left] {} (c10)
    (c30) edge   node[left] {} (c27)
    (c30) edge   node[left] {} (c28)

    ;

	\end{tikzpicture}
	
\end{center}

\begin{landscape} 

$$\begin{array}{|lcl|llll|}
\hline

\multicolumn{7}{|c|}{\textrm{{\bf Table B.3. Degenerations of $4$-dimensional nilpotent commutative  algebras}}}  \\
\hline
\multicolumn{3}{|c|}{\textrm{Degeneration}}  & \multicolumn{4}{|c|}{\textrm{Parametric basis}} \\
\hline
\hline

\mathcal{C}_{02} &\to& \mathcal{C}_{03}  & 
E_{1}^{t}=te_1+e_2,& 
E_{2}^{t}=t^2e_2 + e_3, & 
E_{3}^{t}=t^2e_3, &
E_{4}^{t}=e_4 \\
\hline

\mathcal{C}_{02} &\to& \mathcal{C}_{06}  & 
E_{1}^{t}=te_1,& 
E_{2}^{t}=te_2 -t^{-1}e_3+ t^{-1}e_4, & 
E_{3}^{t}=e_3-e_4, &
E_{4}^{t}=t^2e_4 \\
\hline

\C{03} &\to& \C{07}  & 
E_{1}^{t}=te_1 + \frac{1}{2}t^{-1} e_2,& 
E_{2}^{t}=te_2 - t^{-1} e_4, & 
E_{3}^{t}=e_3 + e_4, &
E_{4}^{t}= -t^2e_4 \\
\hline

\C{04} &\to& \C{01}  & 
E_{1}^{t}=e_1+e_2,& 
E_{2}^{t}=2e_3, & 
E_{3}^{t}=te_2, &
E_{4}^{t}=e_4 \\
\hline

\mathcal{C}_{05} &\to& \mathcal{C}_{02}  & 
E_{1}^{t}=te_1,& 
E_{2}^{t}=t^2e_2, & 
E_{3}^{t}=t^4e_3, &
E_{4}^{t}=e_4 \\
\hline


\C{06} &\to& \C{07}  & 
E_{1}^{t}=e_1 + e_2,& 
E_{2}^{t}=te_2, & 
E_{3}^{t}=e_3 + e_4, &
E_{4}^{t}= te_4 \\
\hline

\C{07} &\to& \C{04}  & 
E_{1}^{t}=te_1,& 
E_{2}^{t}=e_2, & 
E_{3}^{t}=te_4, &
E_{4}^{t}=e_3 \\
\hline

\C{08} &\to& \C{04}  & 
E_{1}^{t}=e_2,& 
E_{2}^{t}=e_3, & 
E_{3}^{t}=e_4, &
E_{4}^{t}=te_1 \\
\hline

\C{09} &\to& \C{11}  & 
E_{1}^{t}=te_1-2e_2-2t^{-2}e_3,& 
E_{2}^{t}=e_3, & 
E_{3}^{t}=te_2-t^{-1} e_3+8t^{-3}e_4, &
E_{4}^{t}=-2e_4 \\
\hline

\C{10} & \to & \C{02}  & 
E_{1}^{t}=e_1,& 
E_{2}^{t}=e_2, & 
E_{3}^{t}=e_4, &
E_{4}^{t}=te_3 \\
\hline

\mathcal{C}_{10} &\to& \mathcal{C}_{11}  & 
E_{1}^{t}=te_1+e_2-\frac{1}{2}t^{-1}e_3,& 
E_{2}^{t}=te_2-e_3, & 
E_{3}^{t}=te_3, &
E_{4}^{t}=t^2e_4 \\
\hline

\C{11} &\to& \C{03}  & 
E_{1}^{t}=te_1,& 
E_{2}^{t}=t^2e_2, & 
E_{3}^{t}=t^2e_3, &
E_{4}^{t}= t^3e_4 \\
\hline

\C{11} &\to& \C{06}  & 
E_{1}^{t}=te_1,& 
E_{2}^{t}=te_2 - e_3, & 
E_{3}^{t}=te_3, &
E_{4}^{t}= e_4 \\
\hline

\C{11} &\to& \C{08}  & 
E_{1}^{t}=-t^3e_2 + t^4e_3,& 
E_{2}^{t}=t^3e_1-\frac 12 te_2+t^2e_3, & 
E_{3}^{t}=t^6e_3, &
E_{4}^{t}=t^8e_4 \\
\hline

\C{12} &\to& \C{09}  & 
E_{1}^{t}=te_1,& 
E_{2}^{t}=it^2e_3, & 
E_{3}^{t}=e_2-ie_3, &
E_{4}^{t}=t^2e_4 \\
\hline

\mathcal{C}_{12} &\to& \mathcal{C}_{10}  & 
E_{1}^{t}=te_1+te_3,& 
E_{2}^{t}=t^2e_2+t^2e_4, & 
E_{3}^{t}=t^3e_3, &
E_{4}^{t}=t^4e_4 \\
\hline

\mathcal{C}_{13} &\to& \mathcal{C}_{09}  & 
 E_{1}^{t}=e_2,& 
 E_{2}^{t}=e_3, & 
 E_{3}^{t}=-te_1, &
 E_{4}^{t}=te_4 \\
\hline

\mathcal{C}_{13} &\to& \mathcal{C}_{10}  & 
E_{1}^{t}=te_1,& 
E_{2}^{t}=t^2e_2, & 
E_{3}^{t}=e_3- t^{-3}e_4, &
E_{4}^{t}=te_4 \\
\hline

\C{13} &\to& \C{20}(-1)  & 
E_{1}^{t}=e_1+t^{-1}e_2,& 
E_{2}^{t}=e_2+t^{-2}e_3, & 
E_{3}^{t}=t^{-1}e_3+t^{-2}e_4, &
E_{4}^{t}=t^{-1}e_4 \\
\hline

\mathcal{C}_{14} &\to& \mathcal{C}_{05}  & 
E_{1}^{t}=te_1,& 
E_{2}^{t}=t^2e_2, & 
E_{3}^{t}=t^4e_3, &
E_{4}^{t}=t^3e_4 \\
\hline


\mathcal{C}_{14} &\to& \mathcal{C}_{13}  & 
E_{1}^{t}=t^{-1}e_1,& 
E_{2}^{t}=t^{-2}e_2, & 
E_{3}^{t}=t^{-4}e_3, &
E_{4}^{t}=t^{-5}e_4 \\
\hline


\mathcal{C}_{14} &\to& \mathcal{C}_{20}(\alpha\not=-1)  & 
 E_{1}^{t}=te_1+\frac{i}{\sqrt{\alpha+1}}te_2,& 
 E_{2}^{t}=t^2e_2-\frac{1}{\alpha+1}t^2e_3+\frac{2i}{\sqrt{\alpha+1}}t^2e_4, & 
 E_{3}^{t}=\frac{i}{\sqrt{\alpha+1}}t^3e_3+\frac{\alpha}{\alpha+1}t^3e_4, &
 E_{4}^{t}=\frac{i}{\sqrt{\alpha+1}}t^4e_4 \\
\hline

\C{14} &\to& \C{25}  & 
E_{1}^{t}=te_1+ite_2-ite_3,& 
E_{2}^{t}=te_1-ite_2+ite_3, & 
E_{3}^{t}=2t^2e_3, &
E_{4}^{t}=2t^3e_4 \\
\hline


\mathcal{C}_{15} &\to& \mathcal{C}_{05}  & 
 E_{1}^{t}=t(e_1+e_3),& 
 E_{2}^{t}=t^2e_2, & 
 E_{3}^{t}=t^4e_3, &
 E_{4}^{t}=t^3e_4 \\
\hline

\mathcal{C}_{15} &\to& \mathcal{C}_{12}  & 
E_{1}^{t}=te_1,& 
E_{2}^{t}=t^2e_2, & 
E_{3}^{t}=i(t^2e_2-te_3+e_4), &
E_{4}^{t}=t^3e_4 \\
\hline

\mathcal{C}_{15} &\to& \mathcal{C}_{20}(\alpha)  & 
 E_{1}^{t}=te_1+te_2+(1-\alpha)te_3,& 
 E_{2}^{t}=t^2e_2+t^2e_3+(2-2\alpha)t^2e_4, & 
 E_{3}^{t}=t^3e_3+(2-\alpha)t^3e_4, &
 E_{4}^{t}=t^4e_4 \\
\hline

\mathcal{C}_{15} &\to& \mathcal{C}_{25}  & 
E_{1}^{t}=te_1+te_2+te_3,& 
E_{2}^{t}=-te_1+te_2+te_3+4te_4, & 
E_{3}^{t}=2t^2e_3+4t^2e_4, &
E_{4}^{t}=2t^3e_4 \\
\hline

\C{16} &\to& \C{14}  & 
 E_{1}^{t}=te_1+t^3e_3,& 
 E_{2}^{t}=t^{2}e_2+2t^4e_4, & 
 E_{3}^{t}=t^{4}e_3, &
 E_{4}^{t}=t^{5}e_4 \\
\hline

\mathcal{C}_{16} &\to& \mathcal{C}_{15}  & 
E_{1}^{t}=t^{-1}e_1,& 
E_{2}^{t}=t^{-2}e_2, & 
E_{3}^{t}=t^{-4}e_3, &
E_{4}^{t}=t^{-6}e_4 \\

\hline

\mathcal{C}_{16} &\to& \mathcal{C}_{22}  & 
E_{1}^{t}=te_1+t(t-1)e_2+t(t-1)^2(t-2)e_3,& 
E_{2}^{t}=t^{2}e_2+t^2(t-1)^2e_3+2t^3(t-1)^2(t-2)e_4, & 
E_{3}^{t}=-t^{3}e_3-2t^3(t-1)^2e_4, &
E_{4}^{t}=-t^{5}e_4 \\
\hline

\C{17} & \to &  \C{05}  & 
 E_{1}^{t}=te_1+te_2,& 
 E_{2}^{t}=t^2e_2+t^2e_3, & 
 E_{3}^{t}=t^{4}e_4, &
 E_{4}^{t}=t^3e_3 \\
\hline

\C{17} & \to &  \C{12}  & 
E_{1}^{t}=te_1,& 
E_{2}^{t}=t^2e_2, & 
E_{3}^{t}=e_3-t^{-4}e_4, &
E_{4}^{t}=e_4 \\
\hline

\C{17} &\to& \C{13}  & 
E_{1}^{t}=te_1+e_3,& 
E_{2}^{t}=t^2e_2+e_4, & 
E_{3}^{t}=t^4e_3, &
E_{4}^{t}= t^4 e_4 \\
\hline

\C{17} &\to& \C{20}(\alpha\not=-1)  & 
E_{1}^{t}=te_1+te_2+\frac{1}{\alpha+1}te_3,& 
E_{2}^{t}=t^2e_2+t^2e_3+\frac{1}{(\alpha+1)^2}t^2e_4, & 
E_{3}^{t}=t^3e_3+\frac{1}{\alpha+1}t^3e_4, &
E_{4}^{t}= \frac{1}{\alpha+1}t^4e_4 \\
\hline

\C{17} &\to& \C{25}  & 
E_{1}^{t}=\frac{1}{2}te_1+\frac{1}{2}te_2+\frac{1}{2}te_3+\frac{1}{2}te_4,& 
E_{2}^{t}=-\frac{1}{2}te_1+\frac{1}{2}te_2+\frac{1}{2}te_3+\frac{1}{2}te_4, & 
E_{3}^{t}=\frac{1}{2}t^2e_3+\frac{1}{2}t^2e_4, &
E_{4}^{t}=\frac{1}{4}t^3e_4 \\
\hline

\C{18} &\to& \C{14}  & 
 E_{1}^{t}=te_1+t^{-1}e_3,& 
 E_{2}^{t}=t^2e_2+t^{-2}e_4, & 
 E_{3}^{t}=t^4e_3, &
 E_{4}^{t}= t^3 e_4 \\
\hline

\C{18} &\to& \C{17}  & 
E_{1}^{t}=t^{-1}e_1+t^{-1}e_2+t^{-2}e_3,& 
E_{2}^{t}=t^{-2}e_2+t^{-2}e_3+(2t^{-2}+t^{-4})e_4, & 
E_{3}^{t}=t^{-4}e_3+t^{-4}e_4, &
E_{4}^{t}= t^{-8} e_4 \\
\hline

\end{array}$$ 

$$\begin{array}{|lcl|llll|}
\hline



\C{18} &\to& \C{22}  & 
 E_{1}^{t}=te_1+t\sqrt[5]{t+1}e_2+\frac{t^2}{\sqrt[5]{(t+1)^2}}e_3,& 
 E_{2}^{t}=t^{2}e_2+t^{2}\sqrt[5]{(t+1)^2}e_3+\frac{2t^2+2t^3+t^4}{\sqrt[5]{(t+1)^4}}e_4, & 
 E_{3}^{t}=t^{3}e_3+t^{3}\sqrt[5]{(t+1)^4}e_4, &
 E_{4}^{t}= t^{5} e_4 \\
\hline

\mathcal{C}_{19}(\af) &\to& \mathcal{C}_{14}  & 
 E_{1}^{t}=e_1+te_2+t^{-1}e_3,& 
 E_{2}^{t}=e_2+(2\alpha t^{-1}+t^{-2})e_4, &
 E_{3}^{t}=e_3, & 
 E_{4}^{t}=t^{-1} e_4 \\
\hline

\mathcal{C}_{19}(0) &\to& \mathcal{C}_{15}  & 
  E_{1}^{t}=te_1,& 
  E_{2}^{t}=t^2e_2, &
  E_{3}^{t}=t^4e_3, & 
  E_{4}^{t}=t^6 e_4 \\
\hline

\mathcal{C}_{19}(\af) &\to& \mathcal{C}_{17}  & 
 E_{1}^{t}=t^{-1}e_1,& 
 E_{2}^{t}=t^{-2}e_2, &
 E_{3}^{t}=t^{-4}e_3, & 
 E_{4}^{t}=t^{-8} e_4 \\
\hline

\mathcal{C}_{19}(\af\neq 0) &\to& \mathcal{C}_{22}  & 
 \begin{array}{l}
\scriptstyle
E_{1}^{t}=tf_1-(\gamma+1)\left(t+\frac{\gamma-1}{5\gamma^2-1}t^2\right)f_2\\
\scriptstyle+(\gamma^2-1)t^2f_3
\end{array},& 

 \begin{array}{l}
\scriptstyle
E_{2}^{t}=-t^2f_2+(\gamma^2-1)\left(t^2+\frac{2\gamma}{5\gamma^2-1}t^3\right)f_3\\
\scriptstyle
+(\gm^2 - 1)^2t^2\left(\gm(\gm + 1)^2+\frac{2\gm(\gm + 1)^2(\gm - 1)}{(5\gm^2 - 1)}t\right.\\
\scriptstyle
\left.+\frac{8\gm^2(3\gm^2 - 1)}{(5\gm^2 - 1)^2}t^2-\frac{4\gm(\gm^2 - 1)^2}{(5\gm^2 - 1)^3}t^3\right.\\
\scriptstyle\left.-\frac{2(3\gm^2 + 1)(\gm^2 - 1)^2}{(5\gm^2 - 1)^4}t^4\right)f_4
\end{array}, &
 
\begin{array}{l}
\scriptstyle
E_{3}^{t}=\gamma t^3f_3+\gamma(\gamma^2-1)^2\\
\scriptstyle
\cdot\left((\gamma+1)t^3+\frac{4\gamma}{5\gamma^2-1}t^4+\frac{2(3\gamma^2+1)}{(5\gamma^2-1)^2}t^5\right)f_4
\end{array}, & 

\scriptstyle
E_{4}^{t}=\gamma(\gamma^2-1) t^5 f_4 \\
&&&&&&\\
&&&\multicolumn{4}{|c|}{\scriptstyle  f_1=ie_1+e_2-(1+\af i)e_3-2(\af+i)^2e_4,\ f_2=e_2-e_3-(\af-i)^2e_4,\ f_3=e_3-(1+\af i)e_4,\ f_4=e_4,\ (\gamma^2-1)^2\gamma=\alpha i,\ 5\gamma^2\not=1}\\
\hline

\mathcal{C}_{19}(0) &\to& \mathcal{C}_{23}  & 
  E_{1}^{t}=te_1+ite_2-ite_3,& 
  E_{2}^{t}=t^2e_2-t^2e_3+t^2e_4, &
  E_{3}^{t}=it^3e_3-it^3e_4, & 
  E_{4}^{t}=-t^6 e_4 \\
\hline

\mathcal{C}_{19}(0) &\to& \mathcal{C}_{29}  & 
E_{1}^{t}=\frac{i}{2}te_1+\frac{1}{2}te_2-\frac{1}{2}te_3,& 
E_{2}^{t}=-\frac{i}{2}te_1+\frac{1}{2}te_2-\frac{1}{2}te_3+(t+t^3)e_4, &
E_{3}^{t}=\frac{1}{2}t^2e_3-\frac{1}{4}(2t^2+t^4)e_4, & 
E_{4}^{t}=\frac{1}{4}t^4 e_4 \\
\hline

\mathcal{C}_{20}(\af\neq 0) &\to& \mathcal{C}_{11}  & 
E_{1}^{t}=te_1+\frac \af 2 e_2-\frac{\af^3}{8}t^{-1}e_3,& 
E_{2}^{t}=\af te_3+\frac{\af^3}2e_4, &
E_{3}^{t}=te_2 - \frac{\af^3}{2}t^{-1}e_4, & 
E_{4}^{t}=\af t^2 e_4 \\
\hline


\C{20}(0) &\to& \C{26}  & 
E_{1}^{t}=te_1,& 
E_{2}^{t}=e_2, &
E_{3}^{t}=te_3, & 
E_{4}^{t}= t^2e_4 \\
\hline

\mathcal{C}_{21}  &\to& \mathcal{C}_{09}  & 
E_{1}^{t}=te_1,& 
E_{2}^{t}=t^2e_2, &
E_{3}^{t}= e_3, & 
E_{4}^{t}=t^2e_4 \\ 

\hline

\C{21} &\to& \C{20}(4)  & 
E_{1}^{t}=e_1+t^{-1}e_2,& 
E_{2}^{t}=e_2+2t^{-1}e_3, &
E_{3}^{t}=e_3+2t^{-2}e_4, & 
E_{4}^{t}= t^{-1}e_4 \\
\hline

\C{21} &\to& \C{27}  & 
E_{1}^{t}=te_2-t^2e_4,& 
E_{2}^{t}=te_1, &
E_{3}^{t}=t^2e_3, & 
E_{4}^{t}= t^3e_4 \\
\hline

\C{22} &\to& \C{05}  & 
 E_{1}^{t}=te_1-te_2,& 
 E_{2}^{t}=t^{2}e_2-2t^{2}e_3, & 
 E_{3}^{t}=-4t^{4}e_4, &
 E_{4}^{t}=t^{3}e_3 \\
\hline
    
\C{22} &\to& \C{12}  & 
E_{1}^{t}=te_1-te_2,& 
E_{2}^{t}=t^{2}e_2-2t^2e_3, & 
E_{3}^{t}=it^2e_2+2it^2e_3, &
E_{4}^{t}=-4t^{4}e_4 \\
\hline

\C{22} &\to& \C{20}(\alpha\not=4)  & 
 E_{1}^{t}=te_1+\frac{\alpha}{4-\alpha}te_2,& 
 E_{2}^{t}=t^{2}e_2+\frac{2\alpha}{4-\alpha}t^{2}e_3, & 
 E_{3}^{t}=t^3e_3+\frac{8\alpha}{(4-\alpha)^2}t^{3}e_4, &
 E_{4}^{t}=\frac{4}{4-\alpha}t^{4}e_4 \\
\hline

\C{22} &\to& \C{21}  & 
E_{1}^{t}=t^{-1}e_1,& 
E_{2}^{t}=t^{-2}e_2, & 
E_{3}^{t}=t^{-3}e_3, &
E_{4}^{t}=t^{-5}e_4 \\
\hline

\C{22} &\to& \C{25}  & 
E_{1}^{t}=te_1,& 
E_{2}^{t}=te_2, & 
E_{3}^{t}=t^2e_3, &
E_{4}^{t}=t^3e_4 \\
\hline





\C{23} &\to& \C{05}  & 
 E_{1}^{t}=te_1+te_2,& 
 E_{2}^{t}=t^{2}e_2+2t^{2}e_3, & 
 E_{3}^{t}=4t^{4}e_4, &
 E_{4}^{t}=t^{3}e_3 \\
\hline

\C{23} &\to& \C{20}(\af)  & 
  E_{1}^{t}=te_1+\frac{\sqrt{\alpha}}{2}te_2+te_3,& 
  E_{2}^{t}=t^2e_2+\sqrt{\alpha}t^2e_3+t^2e_4, & 
  E_{3}^{t}=t^3e_3+\sqrt{\alpha}t^3e_4, &
  E_{4}^{t}=t^{4}e_4 \\
\hline

\C{23} &\to& \C{30}  & 
E_{1}^{t}=te_1,& 
E_{2}^{t}=e_2-e_4, & 
E_{3}^{t}=te_3, &
E_{4}^{t}=t^{2}e_4 \\
\hline


\C{24} &\to& \C{22}  & 
 E_{1}^{t}=te_1+\frac{i}{2}te_2+it^2e_3& 
 E_{2}^{t}=t^2e_2+it^2e_3-(\frac 14 t^2+t^4)e_4, & 
 E_{3}^{t}=t^3e_3+(\frac i2 t^3-t^4)e_4, &
 E_{4}^{t}=it^5e_4 \\
\hline

\mathcal{C}_{24}  &\to& \mathcal{C}_{23}  & 
E_{1}^{t}=t^{-1}e_1,& 
E_{2}^{t}=t^{-2}e_2, &
E_{3}^{t}=t^{-3}e_3, & 
E_{4}^{t}=t^{-6}e_4 \\ 
\hline

\C{24} &\to& \C{29}  & 
E_{1}^{t}=e_1& 
E_{2}^{t}=t^{-1}e_2-t^{-3}e_4, & 
E_{3}^{t}=t^{-1}e_3, &
E_{4}^{t}=t^{-2}e_4 \\
\hline


\mathcal{C}_{25} & \to& \C{09}  & 
E_{1}^{t}=e_1-e_2,& 
E_{2}^{t}=-2e_3, & 
E_{3}^{t}=te_2, &
E_{4}^{t}=-2te_4 \\

\hline

\mathcal{C}_{25} & \to& \C{27}  & 
E_{1}^{t}=2te_1,& 
E_{2}^{t}=2t^2e_2+2t^2e_3, & 
E_{3}^{t}=4t^3e_3+4t^3e_4, &
E_{4}^{t}=8t^4e_4 \\

\hline

\C{26} &\to& \C{03}  & 
E_{1}^{t}=e_1+e_2,& 
E_{2}^{t}=2e_3, & 
E_{3}^{t}=2e_4, &
E_{4}^{t}=te_2 \\
\hline

\mathcal{C}_{27} & \to & \mathcal{C}_{11}  & 
E_{1}^{t}=t^2e_1+2t^2e_2+2t^2e_3,& 
E_{2}^{t}=4t^4e_3+8t^4e_4, & 
E_{3}^{t}=2t^3e_2, &
E_{4}^{t}=4t^6e_4 \\

\hline

\mathcal{C}_{27} & \to & \mathcal{C}_{26}  & 
E_{1}^{t}=e_1,& 
E_{2}^{t}=te_2, & 
E_{3}^{t}=te_3, &
E_{4}^{t}=te_4 \\

\hline

\mathcal{C}_{28} &\to& \mathcal{C}_{02}  & 
E_{1}^{t}=e_1+e_2,& 
E_{2}^{t}=2e_3, & 
E_{3}^{t}=4e_4, &
E_{4}^{t}=te_2 \\
\hline

\C{28} &\to& \C{11}  & 
E_{1}^{t}=t^2e_1+2t^2e_2+t^2e_3,& 
E_{2}^{t}=4t^4e_3+t^4e_4, & 
E_{3}^{t}=8t^3e_2+2t^3e_3, &
E_{4}^{t}=4t^6e_4 \\

\hline


\mathcal{C}_{29} &\to& \mathcal{C}_{12}  & 
E_{1}^{t}=te_1+te_2,& 
E_{2}^{t}=2t^2e_3+2t^2e_4, & 
E_{3}^{t}=2t^2e_2, &
E_{4}^{t}=4t^4e_4 \\
\hline

\mathcal{C}_{29} &\to& \mathcal{C}_{25}  & 
E_{1}^{t}=te_1+ite_3,& 
E_{2}^{t}=te_2+ite_3, & 
E_{3}^{t}=t^2e_3-t^2e_4, &
E_{4}^{t}=it^3e_4 \\
\hline


\mathcal{C}_{29} &\to& \mathcal{C}_{30}  & 
E_{1}^{t}=t^{-1}e_1,& 
E_{2}^{t}=e_2, & 
E_{3}^{t}=t^{-1}e_3, &
E_{4}^{t}=t^{-2}e_4 \\
\hline

\mathcal{C}_{30} &\to& \mathcal{C}_{10}  & 
E_{1}^{t}=te_1+te_2,& 
E_{2}^{t}=2t^2e_3+t^2e_4, & 
E_{3}^{t}=-8t^3e_2, &
E_{4}^{t}=4t^4e_4 \\
\hline


\mathcal{C}_{30} &\to& \mathcal{C}_{27}  & 
E_{1}^{t}=e_1+ie_3,& 
E_{2}^{t}=-it^2e_2+te_3, & 
E_{3}^{t}=-it^{2}e_3+ite_4, &
E_{4}^{t}=t^{2}e_4 \\
\hline

\mathcal{C}_{30} &\to& \mathcal{C}_{28}  & 
E_{1}^{t}=e_1,& 
E_{2}^{t}=t^{-1}e_2, & 
E_{3}^{t}=t^{-1}e_3, &
E_{4}^{t}=t^{-2}e_4 \\
\hline
\end{array}$$

$$\begin{array}{|lcl|llll|}
\hline

\multicolumn{7}{|c|}{\textrm{{\bf Table B.4. Orbit closures  of $4$-dimensional nilpotent commutative  algebras}}}  \\
\hline
\multicolumn{3}{|c|}{\textrm{Degeneration}}  & \multicolumn{4}{|c|}{\textrm{Parametric basis}} \\
\hline
\hline

\mathcal{C}_{19}(t) &\to& \mathcal{C}_{16}  & 
E_{1}^{t}=te_1+t^3e_2+t^5e_3,& 
E_{2}^{t}=t^2e_2+t^6e_3+(2t^7+2t^8+t^{10})e_4, &
E_{3}^{t}=t^4e_3+(2t^8+t^{12})e_4, & 
E_{4}^{t}=t^6 e_4 \\
\hline

\mathcal{C}_{19}(-t^{-5}) &\to& \mathcal{C}_{18}  & 
 E_{1}^{t}=t^{-1}e_1+t^{-6}e_3,& 
 E_{2}^{t}=t^{-2}e_2-t^{-2}e_3+(t^{-2}-t^{-12})e_4, &
 E_{3}^{t}=t^{-4} e_3- t^{-4} e_4, & 
 E_{4}^{t}=t^{-8} e_4 \\
\hline

\mathcal{C}_{19}(-it^2) &\to& \mathcal{C}_{24}  & 
 E_{1}^{t}=te_1+ite_2+i(t^3-t)e_3,& 
 E_{2}^{t}=t^2e_2-t^2e_3+(t^6 - 2t^4 + t^2)e_4, &
 E_{3}^{t}=it^3e_3+i(t^5-t^3)e_4, & 
 E_{4}^{t}=-t^6 e_4 \\
\hline

\mathcal{C}_{20}(t^{-4}) &\to& \mathcal{C}_{05}  & 
E_{1}^{t}=te_1,& 
E_{2}^{t}=t^2e_2, &
E_{3}^{t}= e_4, & 
E_{4}^{t}=t^3e_3 \\ 
\hline

\mathcal{C}_{20}(t^{-1}) &\to& \mathcal{C}_{10}  & 
E_{1}^{t}=te_1,& 
E_{2}^{t}=t^2e_2, &
 E_{3}^{t}= t^2e_3, & 
 E_{4}^{t}=t^3e_4 \\ 
\hline

\mathcal{C}_{20}(t) &\to& \mathcal{C}_{27}  & 
E_{1}^{t}=te_1,& 
E_{2}^{t}=te_2, &
E_{3}^{t}=t^2e_3, & 
E_{4}^{t}=t^3e_4 \\ 
\hline

\end{array}$$

$$\begin{array}{|rcl|l|}
\hline
\multicolumn{4}{|c|}{\textrm{{\bf Table B.5. 
Non-degenerations of $4$-dimensional nilpotent commutative  algebras}}}  \\
\hline
\hline
\multicolumn{3}{|c|}{\textrm{Non-degeneration}} & \multicolumn{1}{|c|}{\textrm{Arguments}}\\
\hline
\hline


 \C{05}  & \not \to&  \C{08},  \C{26}  & {\mathcal R}=
\left\{
\begin{array}{l} A_1A_3=0 \end{array}  \right\} 
 \\
\hline



\C{12} & \not \to & 
\C{26}  & {\mathcal R}=
\left\{ \begin{array}{l} A_1A_4=0, A_1A_2 \subseteq A_4 \\  \end{array}
\right\}
 \\
\hline 

\C{13}  &\not \to& \C{05},\C{20}(\af\ne-1), \C{26}& {\mathcal R}=
\left\{
\begin{array}{l}
A_1A_4+A_2A_3=0, A_1A_3\subseteq A_4, A_1A_2\subseteq A_3\\
c_{11}^2(c_{12}^4c_{22}^3-c_{12}^3c_{22}^4)=c_{13}^4((c_{12}^3)^2-c_{11}^3c_{22}^3)
\end{array}
\right\}
 \\
\hline

\C{14}  &\not \to&  \C{12}&
{\mathcal R}= \left\{ \begin{array}{l}
A_1A_4+A_2A_3=0, A_1A_2\subseteq A_3 \end{array} \right\}\\
\hline

\C{15}  &\not \to&  \C{13}&
{\mathcal R}= \left\{ \begin{array}{l}
A_1A_4+A_3A_3=0, A_1A_3\subseteq A_4,A_1A_2\subseteq A_3,
c_{12}^3c_{23}^4=c_{13}^4c_{22}^3 \end{array} \right\}\\
\hline

\C{16}  &\not \to&  \C{17},\C{28} & {\mathcal R}=
\left\{
\begin{array}{l}
A_1A_4+A_3A_3=0, A_1A_3\subseteq A_4 \end{array}
\right\}
 \\
\hline



\C{18}  &\not \to&  \C{15},\C{28} &  {\mathcal R}=
\left\{
\begin{array}{l}
A_1A_4=0, A_1A_3\subseteq A_4, A_1A_2\subseteq A_3, A_1A_1\subseteq A_2,\\
c_{22}^3(c_{11}^2c_{23}^4+c_{11}^3c_{33}^4)=(c_{12}^3)^2c_{33}^4
\end{array}
\right\}
 \\
\hline


\C{19}(\af\not=0)  &\not \to& \C{15}, \C{28}  & {\mathcal R}=
\left\{
\begin{array}{l}
A_1A_4=0, A_1A_3\subseteq A_4, A_1A_2\subseteq A_3, A_1A_1\subseteq A_2,\\
(c_{12}^3(c_{23}^4)^2+c_{12}^4c_{22}^3c_{33}^4-c_{13}^4c_{22}^3c_{23}^4-c_{12}^3c_{22}^4c_{33}^4)^2(c_{11}^2c_{22}^3)^4c_{33}^4
=\alpha^2(c_{11}^3c_{22}^3c_{33}^4+c_{11}^2c_{22}^3c_{23}^4-(c_{12}^3)^2c_{33}^4)^5
\end{array}
\right\}
 \\
\hline

\C{19}(0)  &\not \to& \C{21}  & {\mathcal R}=
\left\{
\begin{array}{l}
A_1A_4=0, A_1A_3\subseteq A_4, A_1A_2\subseteq A_3,
c_{12}^3((c_{23}^4)^2-c_{22}^4c_{33}^4)=c_{22}^3(c_{13}^4c_{23}^4-c_{12}^4c_{33}^4)
\end{array}
\right\}
 \\
\hline

\C{20}(0) & \not \to&  \C{06}, \C{08} & 
\R=\{ \begin{array}{l} A_2A_2=0 \end{array} \} \\
\hline

\C{20}(\af\not=0)   &\not \to&  \C{02},  \C{26} &{\mathcal R}=
\left\{ \begin{array}{l} A_1A_4+A_2A_3=0, A_1A_3+A_2A_2\subseteq A_4, A_1A_2\subseteq A_3, A_1A_1\subseteq A_2,\\
c_{11}^2c_{22}^4=\af c_{12}^3c_{13}^4 \end{array}
\right\}  \\ \hline

\C{21}  & \not \to &  \C{02}, \C{20}(\af\ne 4)  & 
{\mathcal R}= \left\{
\begin{array}{l}
A_1A_4+A_3A_3=0, A_1A_3+A_2A_2\subseteq A_4, A_1A_2\subseteq A_3, A_1A_1\subseteq A_2,\\
2c_{11}^3c_{23}^4=4c_{12}^3c_{13}^4-c_{11}^2c_{22}^4
\end{array}
\right\}
 \\
\hline



\C{23}  &\not \to& \C{09}  &{\mathcal R}=
\left\{
\begin{array}{l}
A_1A_4=0, A_2A_2\subseteq A_4, c_{22}^4c_{33}^4=(c_{23}^4)^2
\end{array}
\right\}
 \\
\hline

\C{24}  &\not \to& \C{13}, \C{15}  & {\mathcal R}=
\left\{
 \begin{array}{l} A_1A_4=0, A_2A_2\subseteq A_4
 \end{array}
\right\}
 \\
\hline

\C{25}  & \not \to & \C{02}
&  {\mathcal R}=
\left\{
 \begin{array}{l} A_1A_1\subseteq A_3, A_3A_3=0
 \end{array}\right\}
 \\
\hline

\C{28}  & \not \to &  \C{26}  &   {\mathcal R}=
\left\{
\begin{array}{l}
A_1A_4=0, A_1A_3+A_2A_2\subseteq A_4, A_1A_1\subseteq A_3\\
(c_{23}^4)^2=c_{22}^4c_{33}^4,\ c_{33}^4(c_{12}^3c_{11}^4-c_{11}^3c_{12}^4)=c_{13}^4(c_{12}^3c_{13}^4-c_{11}^3c_{23}^4)
\end{array}
\right\}

 \\
\hline

\C{29}  &\not \to& \C{05}, \C{20}(\alpha)  & {\mathcal R}=
\left\{
 \begin{array}{l}  A_1A_1\subseteq A_3
 \end{array}\right\}
 \\
\hline

\hline

\end{array}$$

\end{landscape} 
    


\end{document}